w
\RequirePackage{fix-cm}
\documentclass[stropt,envcountsect,envcountsame]{svjour3}          

\smartqed  
\usepackage{graphicx}
\usepackage{graphicx}
\usepackage{amsmath}
\usepackage{amssymb}
\usepackage{amsfonts}
\usepackage{mathtools}
\usepackage[english]{babel}
\usepackage{enumerate}
\usepackage{url}
\usepackage{color}
\usepackage{subfigure}
\usepackage[normalem]{ulem}
\usepackage{lipsum}
\usepackage{authblk}
\usepackage{hyperref}
\usepackage{cite}
\usepackage{todonotes}
\usepackage{mathabx}
\usepackage{tikz}
\usepackage{bm}
\usepackage{pgfplots}

\DeclareMathOperator    \cone            {cone}

\DeclareMathOperator    \intr               {int}

\DeclareMathOperator    \rec                {rec}

\newcommand{\R}{\mathbb R}
\newcommand{\Q}{\mathbb Q}
\newcommand{\Z}{\mathbb Z}
\newcommand{\N}{\mathbb N}

\def\ve#1{\mathchoice{\mbox{\boldmath$\displaystyle\bf#1$}}
{\mbox{\boldmath$\textstyle\bf#1$}}
{\mbox{\boldmath$\scriptstyle\bf#1$}}
{\mbox{\boldmath$\scriptscriptstyle\bf#1$}}}

\newcommand{\x}{{\ve x}}

\newcommand{\ceil}[1]{\lceil #1 \rceil}

\newcommand{\normone}[1]{\left\lVert#1\right\rVert_1}
\newcommand{\norm}[1]{\left\lVert#1\right\rVert}

\newcommand{\C}{\mathcal C}

\newtheorem{observation}[theorem]{Observation} 

\def\journalversion{Include-Proofs}

\begin{document}
\title{Complexity, Exactness, and Rationality in Polynomial Optimization \thanks{A.~Del Pia is partially funded by ONR grant N00014-19-1-2322. D. Bienstock is partially funded by ONR grant N00014-16-1-2889. R. Hildebrand is partially funded by ONR Grant N00014-20-1-2156 and by AFOSR grant FA9550-21-0107.  Any opinions, findings, and conclusions or recommendations expressed in this material are those of the authors and do not necessarily reflect the views of the Office of Naval Research or the Air Force Office of Scientific Research.}}
\titlerunning{Exact polynomial optimization}
%

\author{Daniel Bienstock
\and
Alberto Del Pia
\and
Robert Hildebrand
}
\authorrunning{Bienstock et al.}
%
\institute{
Daniel Bienstock, Department of Industrial Engineering and Operations Research, Columbia University 
\email{dano@columbia.edu}
\and
Alberto Del Pia, Department of Industrial and Systems Engineering \& Wisconsin Institute for Discovery,
             University of Wisconsin-Madison 
\email{delpia@wisc.edu }
\and 
Robert Hildebrand, Grado Department of Industrial and Systems Engineering,
Virginia Tech 
\email{rhil@vt.edu}
}
\maketitle              
\begin{abstract}
We focus on rational solutions or nearly-feasible rational solutions that serve as certificates of feasibility for polynomial optimization problems.  We show that, under some separability conditions, certain cubic polynomially constrained sets admit rational solutions.  However,  
 we show in other cases that it is NP Hard to detect if rational solutions exist or if they exist of any  reasonable size.  We extend this idea to various settings including near feasible, but super optimal solutions and detecting rational rays on which a cubic function is unbounded. Lastly, we show that in fixed dimension, the feasibility problem over a set defined by polynomial inequalities is in NP by providing a simple certificate to verify feasibility.  We conclude with several related examples of irrationality and encoding size issues in QCQPs and SOCPs.
 \keywords{Polynomial Optimization \and Algebraic Optimization \and Rational solutions \and NP \and Cubic \and Quadratic}
\end{abstract}

\section{Introduction}
This paper addresses basic questions of precise certification of feasibility and optimality, for optimization problems with polynomial constraints, in polynomial time, under the Turing model of computation.  Recent progress in polynomial optimization and mixed-integer nonlinear programming has produced elegant methodologies and effective   implementations; however such implementations may produce \textit{imprecise} solutions whose actual quality can be difficult to rigorously certify, even approximately. The work we address is motivated by these issues, and can be summarized as follows:

\noindent {\bf Question:} given a polynomially constrained problem, what can be said about the existence of feasible or approximately feasible rational solutions of polynomial size (bit encoding length)\footnote{
Throughout we will use the concept of \emph{size}, or bit encoding length, of rational numbers,  vectors, linear inequalities, and formulations. For these standard definitions we refer the reader to Section 2.1 in \cite{Schrijver}}, and more generally the existence of rational, feasible or approximately-feasible solutions that are also approximately-optimal for a given polynomial objective? 

As is well-known, Linear Programming is polynomially solvable \cite{KHACHIYAN198053, Karmarkar1984}, and, moreover, every face of a rational polyhedron contains a point of polynomial size
~\cite{Schrijver}. If we instead optimize a quadratic function over linear constraints, the problem becomes NP-Hard~\cite{Pardalos1991}, but perhaps surprisingly, Vavasis~\cite{Vav90} proved that a feasible system consisting of linear inequalities and just one quadratic inequality, all with rational coefficients, always has a rational feasible solution of polynomial size. This was extended by Del Pia, Dey, and Molinaro~\cite{Pia2016} to show that the same result holds in the mixed-integer setting.  See also \cite{Hochbaum2007} discussion of mixed-integer nonlinear optimization problems with linear constraints. 

On the negative side, there are classical examples of SOCPs all of whose feasible solutions require exponential size \cite{alizadeh}, \cite{ramana}, \cite{letchfordparkes}.  Adding to this, in Example~\ref{ex:irrational}, we provide an SOCP where all feasible solutions are irrational (likely, a folklore result). Recent results of Pataki and Touzov~\cite{Pataki2021} actually show that many SDP's have large encoding size issues.  Their work hinges on earlier examples from Khachiyan.  In the nonconvex setting, there are examples of quadratically constrained, linear objective problems, on $n$ bounded variables, and with coefficients of magnitude $O(1)$, that admit solutions with maximum additive infeasibility $O(2^{-2^{\Theta(n)}})$ but multiplicative (or additive) superoptimality $\Theta(1)$. O'Donnell \cite{odonnell} questions whether SDPs associated with fixed-rank iterates of sums-of-squares hierarchies (which relax nonconvex polynomially constrained problem) can be solved in polynomial time, because optimization certificates might require exponential size.  The  issue of accuracy in solutions is not just of theoretical interest. As an example,  \cite{wakistrange} 
describes instances of SDPs (again, in the sums-of-squares setting)  where a solution is very nearly certified as optimal, and yet proves substantially suboptimal.
 
Vavasis' result suggests looking at systems of two or more quadratic constraints, or (to some extent equivalently) optimization problems where the objective is quadratic, and at least one constraint is quadratic, with all other constraints linear. The problem of optimizing a quadratic subject to one quadratic constraint (and no linear constraints) can be solved in polynomial time using semidefinite-programming techniques \cite{polter}, to positive tolerance. When the constraint is positive definite (i.e. a ball constraint) the problem can be solved to tolerance $\epsilon$ in time $\log \log \epsilon^{-1}$ \cite{advances}, \cite{karmarkar} (in other words
$O(k)$ computations guarantee accuracy $2^{-2^k}$).  Vavasis \cite{vavasiszippel} proved, on the other hand, that \textit{exact} feasibility of a system of two quadratics can be tested in polynomial time.  

With regards to systems of more than two quadratic constraints, Barvinok \cite{Barvinok} proved a fundamental result: for each fixed integer $m$ there is an algorithm that, given $n \times n$ rational matrices  $A_i$ ($1 \le i \le m$) tests, in polynomial time,
feasibility of the system of equations 
$$ x^T A_i x = 0 \quad \text{for $1 \le i \le m$}, \quad x \in \R^n, \ \|x\|_2 = 1.$$
A feature of this algorithm is that   certification does not
rely on producing a feasible vector; indeed, all feasible solutions may be irrational.  As a corollary of this result, \cite{Bienstock2016} proves that, for
each fixed integer $m$ there is an algorithm that solves, in polynomial time,
an optimization problem of the form
$$ \min f_0(x), \quad \text{s.t.} \quad f_i(x) \le 0 \quad \text{for $1 \le i \le m$} $$
where for $0 \le i \le m$, $f_i(x)$ is an $n$-variate quadratic polynomial, and we assume that 
the quadratic part of $f_1(x)$ positive-definite; moreover a rational vector that is (additively) both $\epsilon$-feasible and -optimal can be computed in time polynomial in the size of the formulation and $\log \epsilon^{-1}$.  An important point with regards to \cite{Barvinok} and \cite{Bienstock2016} is that the analyses do not apply to systems of arbitrarily many \textit{linear} inequalities and just two quadratic inequalities. 

De Loera et. al.~\cite{DELOERA20111260} use the \emph{Nullstellensatz} to provide feasibility and infeasibility certificates to systems of polynomial equations through solving a sequence of large linear equations.  Bounds on the size of the certificates are obtained
\cite{GRIGORIEV2001153}.  This technique does not seem amenable to systems with a large number of linear inequalities due to the necessary transformation into equations and then blow up of the number of variables used. Another approach is to use the Positivestellensatz and compute an infeasibility certificate using sums of squares hierarchies.  As mentioned above, see \cite{odonnell} for a discussion if exactness and size of these hierarchies needed for a certificate.

 Renegar~\cite{Ren92f} shows that the problem of deciding whether a system of polynomial inequalities is feasible can be decided in polynomial time provided that the dimension is considered fixed. This is a landmark result, however, the algorithm and techniques are quite complicated. Renegar~\cite{renegar92} then shows how to provide \emph{approximate solutions} that are near feasible solution. Another technique to obtain a similar result is \emph{Cylindrical Algebraic Decomposition}.  See, e.g.,~\cite{Basu2006}.  These techniques can admit a \emph{rational univariate representation}~\cite{Rouillier1999}, encoding feasible solutions as roots of univariate polynomials.  See~\cite{Basu-Survey-2014} for more recent results and improvements.  In our work, we aim to avoid these techniques and provide an extremely simple certificate that shows the feasibility question is in NP.  

\paragraph{\textbf{Our Results.}}
The main topic we address in this paper is whether a system of polynomial inequalities admits rational, feasible or near-feasible solutions of polynomial size.  We show that it is strongly NP-hard to test if a system of quadratic inequalities that has feasible rational solutions, admits feasible rational solutions of polynomial size (Theorem \ref{hard1}). And it is also hard to test if a feasible system of a linear inequalities and a single cubic inequality has a rational solution (Theorem \ref{hard2}). Next we show that these effects can be seen largely when we consider nearly feasible solutions. We show that a point that is slightly infeasible can have a far superior objective function than any feasible point
(Theorem~\ref{infvsuper}). We then consider unbounded problems. We establish that it is NP-hard to determine whether or not the ray is irrational (Theorem~\ref{hard-unbounded}). In the next section, we show that, given a system of polynomial inequalities on $n$ variables that was known to have a bounded, nonempty feasible region, we can produce as a certificate of feasibility a rational, near-feasible vector that has polynomial size, for fixed $n$ (Theorem \ref{th short certificate}). 
This certificate yields a direct proof that, in fixed dimension, the feasibility problem over a system of polynomial inequalities is in NP.

Next we consider the lower boundaries of Theorem \ref{hard2}. We show that in dimension 2, with one separable cubic inequality and linear inequalities, there exists a rational solution of polynomial size (Theorem~\ref{th onetwo}).  
We then provide background theory maximizing a cubic function over a polyhedron.
We conclude with several related examples of simple sets with complicated solutions. 
 
\section{NP-Hardness of determining existence of rational feasible solutions}
\label{sec:rational-NP-Hard}
 \label{mainnphard}
 In this section we show a number of hardness results concerning systems of polynomial inequalities. 
 In particular we prove that given a 3-SAT formula F there is a polynomial-size system S of polynomial inequalities that always admits a feasible rational solution, and with the property that if F is not satisfiable then every feasible rational solution to S has exponential size, whereas if F is satisfiable then S has a feasible rational solution of linear size (Theorem \ref{hard1}). As a result, colloquially, it is strongly NP-hard to test whether a system of quadratic inequalities which is known to have feasible rational solutions, admits feasible rational solutions of polynomial size.     The same proof technique
shows that it is hard to test whether a feasible system of
quadratic inequalities has a rational solution (Theorem \ref{hard2}).  Finally we show that it is NP-hard to decide if a nearly-feasible solution to a polynomial optimization problem is also 'very' superoptimal -- a precise statement is given in Theorem \ref{infvsuper}.\\

\noindent \textit{Basic notation.}
We denote by $\Z[x_1, \dots, x_n]$ the set of all polynomial functions from $\R^n$ to $\R$ with integer coefficients.
For ease of notation, we write a polynomial $g \in \Z[x_1, \dots, x_n]$ of degree $d$ in the form $g(x) = \sum_{I \in \N^n, \normone I \le d} c_I x^I$, where each $c_I \in \Z$ and $x^I := \prod_{i=1}^n x_i^{I_i}$.\\

\noindent \textit{Two important constructions.}  
Examples \ref{example cubic v2} and \ref{ex:small_solutions} will be used throughout our proofs. They were, to the best of our knowledge, previously unknown.

\begin{figure}
    \centering
a.) \includegraphics[scale=0.23]{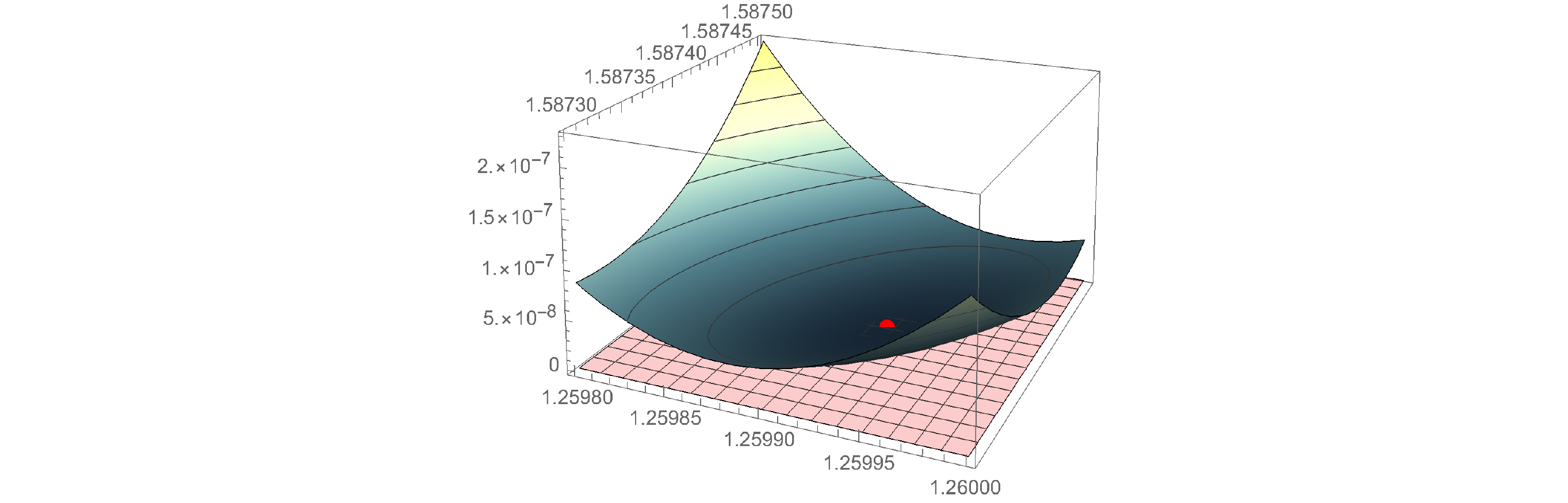} \ 
\includegraphics[scale=0.4]{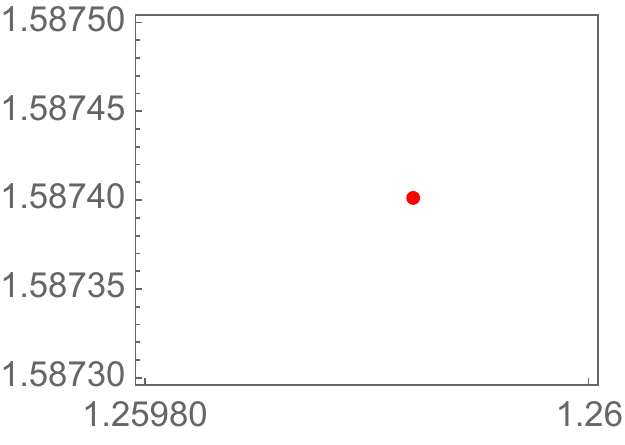}\\
b.) \includegraphics[scale=0.23]{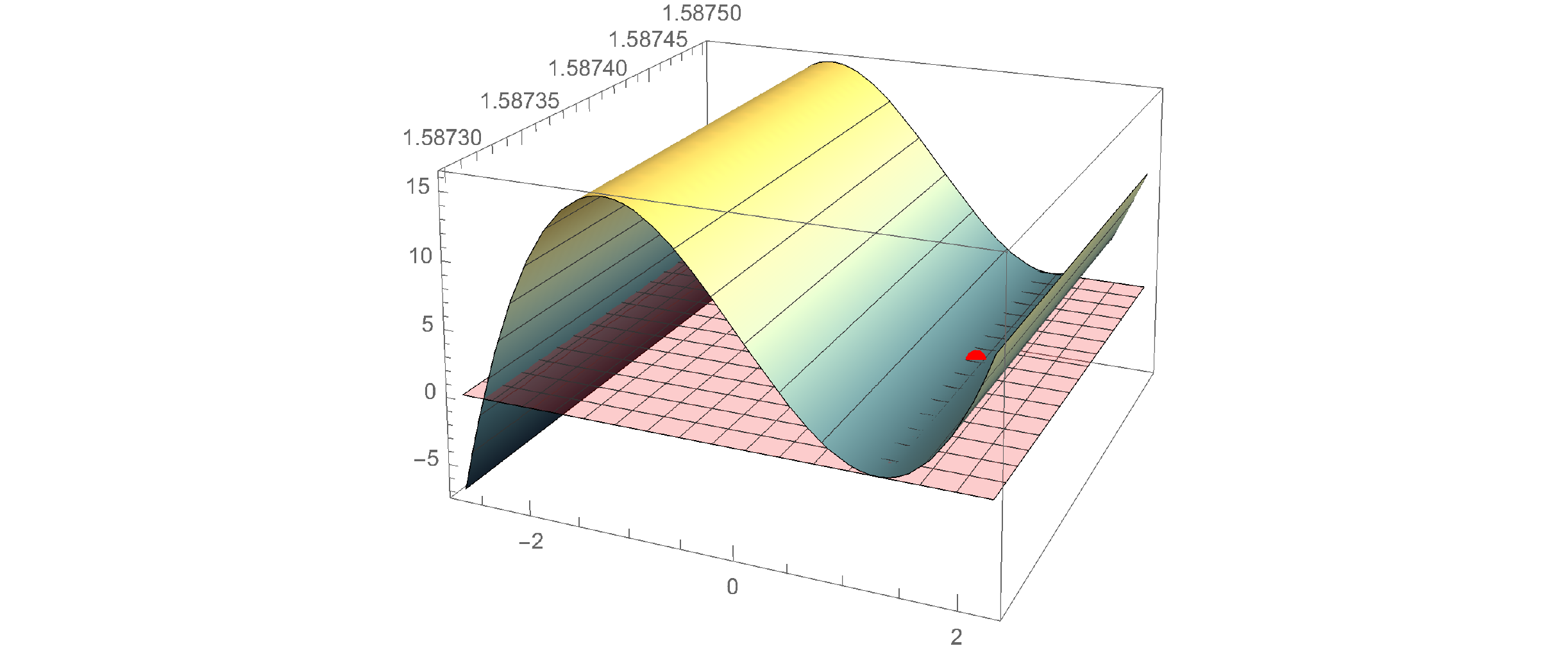} \ \includegraphics[scale=0.4]{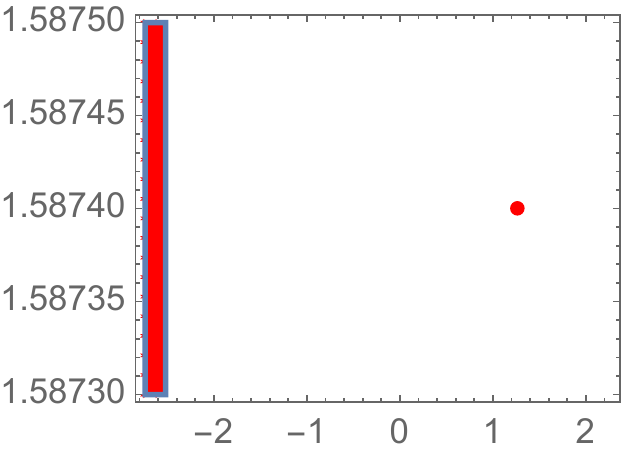}
\caption{The function $h(y)$ on the left on the domain $R_\gamma$, and on the right the set of feasible pints satisfying $y \in R_\gamma$, $h(y) \leq 0$. Row a.) has $\gamma = 0$ and row b.) has $\gamma = 4$.}
    \label{fig: function h}
\end{figure}

\label{sec:special}
\begin{example}[Feasible system with no rational feasible vector]
\label{example cubic v2}
Define   
\begin{equation}
    h(y) := 2y_1^3 + y_2^3 - 6y_1y_2 + 4, \ \ R_\gamma :=  [1.259 - \gamma, 1.26] \times [1.587, 1.59].
\end{equation}
See Figure~\ref{fig: function h} for an illustration.
Then  $\{y \in R_0 : h(y) \leq 0\} = \{y^*\}$ where 
\begin{equation}
y^* := (2^{\frac{1}{3}}, 2^{\frac{2}{3}}) \approx (1.2599, 1.5874). \label{eq:ystardef}
\end{equation} 
\hfill $\diamond$
\end{example}

\begin{observation}
\label{4to12}
\label{obgamma}
\label{ob1}
(1) The point $y^* \in R_\gamma$ for all $\gamma \geq 0$. (2) $y^*$ is the unique minimizer of $h(y)$ for $y \in \R^2_+$.
(3) For $y \in R_4$, $h(y) > - 12$.  (4) The point $\bar y = (-2.74, 1.588) \in R_4$ attains $h(\bar y) < -7$. (5)  $h(y^*) = 0$.
\end{observation}
\begin{proof} (1) is clear. (2) On the curve defined by $y_2 = y_1^2$ --which includes the point $y^*$-- we have that $h(y) = y_1^6 - 4 y_1^3 + 4$ whose sole minimizer is at $y_1 = 2^{1/3} = y^*_1$.  Moreover, $\frac{\partial h}{\partial y_1} = 6(y_1^2 - y_2)$.  These two facts imply that, for any $y \in \R^2_+$, 
$h(y) \ge h(y_2^{1/2}, y_2) \ge h(y^*)$, with at least one of the two inequalities strict unless $y = y^*$as desired. (3) By (2), if $y \in \R^2_+$ then $h(y) \ge h(y^*) = 0$. Otherwise, $y_1 < 0$. Then $\frac{\partial h}{\partial y_1} = 6(y_1^2 - y_2)$ and so either (when $y_1^2 < y_2$) $h(y) \ge h(0, y_2) > 0$, or (when $y_1^2 \ge y_2$) $h(y) \ge h(1.259 - 4, y_2) > -12 + y_2^3  > - 12$. (4), (5) are clear.
\qed
\end{proof}

Note that the proof of Observation~\ref{obgamma} shows that for $R_{0}$, the only zero of $h$ is $y^{*}$, and $h$ takes negative values in $R_{4}$.
In the next example we consider a special system of quadratic inequalities.

\begin{example}\label{ex:small_solutions}
Let $D \subseteq \R \times \R^n$ be the set of $s \in \R$ and $d \in \R^n$ such that the following inequalities are satisfied:
\begin{equation}\label{tiny}
\begin{array}{rlcrl}
 0 & \leq &d_1 &\leq& \frac{1}{2}\\
 0 & \leq& d_{k} &\leq& d_{k-1}^2, \quad k = 2, \ldots, n\\
0 & \leq& s &\leq& d_n^2.
  \end{array}
 \end{equation}
 \end{example}
 
\begin{observation}

\label{ex:small_solutions_obs}

\noindent In any feasible solution to system \eqref{tiny} we either have $s = 0$, or $0 < s \le 2^{-2^n}$ and in this case if $s$ is rational then we need at least $2^n$ bits to represent it using the standard bit encoding scheme. Further, there are rational solutions to system \eqref{tiny} with $s > 0$.\\
We remark that the literature abounds with examples of SOCPs all of whose solutions are doubly exponentially large, see e.g. \cite{alizadeh}, \cite{ramana}.  Our example is similar, however it is \textit{non-convex}. \hfill $\diamond$
\end{observation}

\ifdefined\badjournalversion 
\begin{observation}\label{ob1}
$y^*$ is the unique minimizer of $h(y)$ within $\R^2_+$.
\end{observation}
\begin{proof} On the curve defined by $y_2 = y_1^2$, which includes the point $y^*$, $h(y) = y_1^6 - 4 y_1^3 + 4$ whose sole minimizer is at $y_1 = 2^{1/3} = y^*_1$.  Moreover, $\frac{\partial A}{\partial y_1} = 6(y_1^2 - y_2)$.  These two facts imply that, for any $y \in \R^2_+$, 
$h(y) \ge h(y_2^{1/2}, y_2) \ge h(y^*)$, with at least one of the two inequalities strict unless $y = y^*$as desired. 
\qed
\end{proof}
\fi

\subsection{NP-hardness construction.}
\label{sec:hardintro}
Here we provide our NP-hardness constructions. The main reduction is from the problem 3SAT. An instance of this problem is defined by $n$ \textit{literals} $w_1, \ldots, w_n$ as well as their negations $\bar w_1, \ldots, \bar w_n$, and a set of $m$ clauses
$C_1, \ldots C_m$ where each clause $C_i$ is of the form $(u_{i1} \vee u_{i2} \vee u_{i3})$.  Here, each $u_{ij}$ is a literal or its negation, and $\vee$ means `or'.  The problem is to find `true' or `false' values for each literal, and corresponding values for their negations, so that the formula
\begin{align} \label{eq:3sat}
  & C_1 \wedge C_2 \ldots \wedge C_m
\end{align}
is true, where $\wedge$ means `and'.

 Given an instance of 3SAT as above, we construct a system of quadratic inequalities on the following $3n + 5$ variables:
\begin{itemize}
\item For each literal $w_j$ we have a variable $x_j$; for $\bar w_j$ we use
  variable $x_{n+j}$.
  \item Additional variables $\gamma, \Delta, y_1, y_2, s,$ and $d_1, \ldots, d_n$.
\end{itemize}

We describe the constraints in our quadratically constrained problem\footnote{Constraint \eqref{eq:nasty} as written is cubic, but is equivalent to three quadratic constraints by defining new variables $y_{12} = y_1^2, \quad y_{22} = y_2^2$, and rewriting the constraint as $-n^5 \sum_{j = 1}^{n} x_j^2 \, + \, y_{12} y_1 + y_{22} y_2 - 6 y_1y_2 + 4 \, - \, s \ \le \ -n^6$}. For each clause $C_i = (u_{i_1} \vee u_{i_2} \vee u_{i_3})$ we denote by $x_{i_k}$ the variable associated with the literal $u_{i_k}$ for $1 \le k \le 3$.  In the 
constraints below, the function $h(y)$ and the region $R_\gamma$ are as in Example \ref{example cubic v2}, and
the set $D$ is as in Example \ref{ex:small_solutions}.
\begin{subequations}
\label{eq:np-hard}
\begin{align}
   -1 \le x_j \le 1&& \text{ for } j \in [2n],\label{eq:box}\\
   \quad x_j + x_{n+j} = 0, && \text{for  } j \in [n],\label{eq:box2}\\
   x_{i_1} + x_{i_2} + x_{i_3} \, \ge \, -1 - \Delta, && \text{ for each clause $C_i = (u_{i_1} \vee u_{i_2} \vee u_{i_3})$}\label{eq:join}
\\
   0 \le \gamma, \ 0 \le \Delta \le 2, \ \Delta + \frac{\gamma}{2} \le 2,\label{eq:gammadelta}
\\
          (y_1, y_2) \in R_\gamma, \label{eq:yinRgamma}\\
   (s,d) \in D\label{eq:nastytiny},\\
        -n^5 \textstyle\sum_{j = 1}^{n} x_j^2 \, + \, h(y) - s \leq - n^6. \label{eq:nasty} 
  \end{align}
  \end{subequations}

\noindent As a remark on this system, we note that the feasible region is contained in a bounded region -- in fact, in every feasible solution the absolute value of every variable is at most $4$.  The following observation will be used below.
\begin{observation}\label{lowerboundedx} Suppose $(x, \gamma, \Delta, y, d, s)$ is feasible for \eqref{eq:np-hard}.  Then 
    $$|x_j| \ \geq \ 1- \frac{12}{n^5} - \frac{2^{-2^n}}{n^5} \quad \text{for $1 \leq j \leq 2n$}.$$
\end{observation}
\begin{proof}
 Let $\sigma^2 \coloneq \min_j \{x_j^2\}$.  Since $-1 \leq x_j \leq 1$, we have $1 \geq |x_j| \geq |x_j|^2 \geq \sigma^2$ for all $j = 1, \dots, n$.  Hence by \eqref{eq:nasty}  
 $-n^5(n-1) -n^5 \sigma^2 \ + h(y) -s \le n^6,$ and so
 $$\sigma^2 \ge 1 + \frac{h(y) - s}{n^5}.$$
 Note that $h(y) \geq- 12$ follows from $\gamma\leq 4$ (due to \eqref{eq:gammadelta}) and part 3 of Observation \ref{obgamma} and that\eqref{eq:nastytiny} implies 
  $s\displaystyle \leq 2^{-2^{n}}$. Thus 
  $$\sigma^2 \ge 1 - \frac{12}{n^5} - \frac{2^{-2^n}}{n^5}. \text{\qed}$$
 \qed
\end{proof}

\begin{theorem}\label{hard1} Let $n \ge 3$.  Consider an instance \eqref{eq:3sat}  of 3-SAT and the corresponding  system \eqref{eq:np-hard}.  
\begin{itemize}
    \item [\textbf{(a)}] System \eqref{eq:np-hard} has a rational feasible solution.
    \item [\textbf{(b)}] Suppose formula \eqref{eq:3sat} is satisfiable. Then  \eqref{eq:np-hard} has a rational feasible solution of size at $4n + \mathcal{L}$, where $\mathcal{L}$ is a fixed constant independent of $n$.
    \item [\textbf{(c)}] Suppose formula \eqref{eq:3sat} is not satisfiable. Then  every  rational feasible solution to \eqref{eq:np-hard} has  size at least $2^n$.
\end{itemize}
\end{theorem}

\begin{proof} \hspace{.1in} \\ 
  \noindent \textbf{(a)}  
Set $x^*_j = 1 = -x^*_{n+j}$ for $j \in [n]$, $\Delta^* = 2$, $\gamma^* = 0$,
$d^*_k = 2^{-2^{k-1}}$ for $k \in [n]$ and $s^* = 2^{-2^n}$.  By inspection these
rational values satisfy \eqref{eq:box}, \eqref{eq:box2}, \eqref{eq:join}, \eqref{eq:gammadelta}, and \eqref{eq:nastytiny}.  Further, let $y^* = (2^{\frac{1}{3}}, 2^{\frac{2}{3}})$. Since (Observation \ref{ob1})
\begin{enumerate}
    \item $y^*$ is in the interior of $R_0$,
    \item $y^*$ is the unique local minimizer of $h$,
    \item  $h(y^*) = 0$ and $s^* > 0$, 
    \item $h$ is continuous, and
    \item   the set $\Q^2$ is dense in $\R^2$,
\end{enumerate}
 there exists rational $\hat y \in R_0$ such that $h(\hat y) \leq s^*$.  Hence, \eqref{eq:yinRgamma} and \eqref{eq:nasty} are also satisfied by $(x^*, \gamma^*, \Delta^*, \hat y, d^*, s^*) \in \Q^{3n + 5}$. \\
 
 \noindent \textbf{(b)} Let $\tilde w$ denote a truth assignment that satisfies \eqref{eq:3sat}. 
  For $1 \le j \le n$ set $\tilde x_j = 1 = - \tilde x_{n+j}$ if $\tilde w_j$ is true, else set
  $\tilde x_j = -1 = -\tilde x_{n+j}$.  Set $\tilde \Delta = 0$, $\tilde \gamma = 4$, and $\tilde d_1 = \ldots = \tilde d_n = \tilde s = 0$.  Finally set $(\tilde y_1, \tilde y_2) = (-2.74, 1.588)$. By inspection, and Observation \ref{obgamma}(4) vector $(\tilde x, \tilde \gamma, \tilde \Delta, \tilde y, \tilde d, \tilde s)$ is feasible for \eqref{eq:np-hard}. Its  size is $4n + \mathcal{L}$, where $\mathcal{L}$ is the size of the encoding of the rational numbers $0$, $4$, $0$, $-2.74$ and $1.588$.\\

\noindent \textbf{(c)} 
  Let $(x, \gamma, \Delta, y, d, s)$ be rational feasible. For $1 \le j \le n$ set $w_j$ to be true if $x_j > 0$ and false otherwise.  It follows that there is at least
  one clause $C_i = (u_{i1} \vee u_{i2} \vee u_{i3})$ such that every $u_{ik}$ (for $1 \le k \le 3$) is false, i.e. each $x_{i_k} < 0$. Using constraint \eqref{eq:join} and Observation \ref{lowerboundedx},
  we obtain
  $$ -3 + \frac{36}{n^5} + 3\frac{2^{-2^n}}{n^5} \ge -1 - \Delta, \text{ and by \eqref{eq:gammadelta} } \ \ 
\gamma \le \frac{72}{n^5} + 6\frac{2^{-2^n}}{n^5} \, < \, 1.$$
  This fact has two implications.  First, since $(y_1, y_2) \in R_\gamma \subset \R^2_+$,
  Observation \ref{ob1} implies $h(y) > 0$ (because $y$ is rational). Second,
  constraint \eqref{eq:nasty}, i.e. $-n^5 \sum_{j = 1}^{n} x_j^2 \, +  h(y) \, - \, s \le -n^6$ implies $h(y) \le s$  because $\sum_{j=1}^n x_{j}^{2}\leq\sum_{j=1}^{n}1=n$. So $s > 0$ and since $s \le 2^{-2^n}$
  (by constraint \eqref{eq:nastytiny} or see Observation \eqref{ex:small_solutions_obs}) the proof is complete. 
\qed
\end{proof}
As a summary, we have:
\begin{corollary} \label{cor:firstnphard} System \eqref{eq:np-hard} always has rational feasible solutions. Either it has a rational feasible solution of linear size, or every rational feasible solution has size $\ge 2^n$, and it is strongly NP-Hard to decide which is the case.
\end{corollary}

As an (easy) adaptation of Theorem \ref{hard1} and Corollary~\ref{cor:firstnphard} we have the following theorem.
\begin{theorem}\label{hard2} It is strongly NP-hard to test if there exists a rational solution to a system of the form
$$
f(x) \leq 0 , \quad Ax \leq b,
$$
where $f \in \Z[x_1, \dots, x_n]$ is of degree $3$, and $A \in \Q^{m\times n}$, $b \in \Q^m$.
\end{theorem}
\noindent {\em Proof sketch.} We proceed with a transformation from 3SAT just
as above, except that we dispense with the variables $d_1, \ldots, d_n$ and $s$ and constraint
\eqref{eq:nastytiny} 
and, rather than constraint \eqref{eq:nasty} we impose
  \begin{align}\label{eq:nastytoo}
   & -n^5 \sum_{j = 1}^{n} x_j^2 \, + \, h(y) \ \le \ -n^6.
  \end{align}
  After these changes we are left with a system consisting the linear inequalities \eqref{eq:box}-\eqref{eq:yinRgamma} plus one inequality of degree three, namely \eqref{eq:nastytoo}. Parts (a) and (b) of Theorem \ref{hard1} have identical counterparts: namely, the system of inequalities has a rational solution, and if formula \eqref{eq:3sat} is satisfiable then the system has a rational solution.  Instead of part (c) we argue that if formula \eqref{eq:3sat} is 
  not satisfiable, then in any feasible solution $(x, \gamma, \Delta, y)$ we have $y = y^*$ (and thus, the solution is not rational). To
  do so we proceed as in (c) of Theorem \ref{hard1} to conclude that
   $y \in \R^2_+$ while also $h(y) \le 0$ (no term $-s$ in \eqref{eq:nastytoo}) which yields that $y = y^*$.
\qed   

\subsection{Infeasibility vs superoptimality}

Our next result addresses the interplay between infeasibility and \textit{super}optimality. Consider a polynomial optimization problem
\begin{subequations} \label{eq:genpoly}
\begin{align} 
    \max \ & \ g(x), \quad \text{subject to} \\
     & \ f_i(x) \le 0, \ 1 \le i \le m, \quad x \in \R^n, \label{eq:genpolyconst}
\end{align}
\end{subequations}
where $g$ and the $f_i$ are polynomials. Many popular methods for addressing \eqref{eq:genpoly} focus on numerical algorithms for obtaining (empirically) good solutions, sometimes (often, perhaps) without a quality guarantee; for example IPOPT \cite{IPOPT} and Knitro \cite{Knitro}.  One could also include the method of `rounding' a solution to e.g. a sum-of-squares or Lasserre relaxation, to the nearest rank-one solution. 

Suppose that $\tilde x$ is such a candidate solution to \eqref{eq:genpoly}, and suppose that $g^U$ is a known \textit{upper}-bound\footnote{We stress that \eqref{eq:genpoly} is a maximization problem.} for \eqref{eq:genpoly} such that $g^U - g(\tilde x)$ is small.  In such a case we are likely to characterize $\tilde x$ as `near-optimal' or even `optimal' (if the objective value gap is small enough).  This paradigm is often found in the literature.

However, the nature of the numerical algorithms cited above is such that the vector $\tilde x$ may be slightly infeasible; we further stress that the high-quality algorithms described above often produce very small infeasibilities. Nevertheless, the nonconvex nature of \eqref{eq:genpoly} gives rise to another possibility, namely that the slightly infeasible point $\tilde{x}$ is \textit{super}-optimal, in fact the superoptimality of $\tilde x$ could potentially be very large.

In light of the examples given in previous sections, it is not surprising that cases of \eqref{eq:genpoly} where a vector with very small infeasibility but very large superoptimality can be easily constructed.  Below we prove a stronger result; namely that it is (strongly) NP-hard to decide, given 
a rational vector $\tilde x$ with very small infeasibility, whether the vector attains very large superoptimality or is near-optimal.  Our formal result is given in Theorem \ref{infvsuper} below. We will first make precise the meaning of `small' infeasibilities and `large' superoptimality.

\noindent{\bf Definition.}
Consider an instance of problem \eqref{eq:genpoly} and let $x \in \R^n$. 
Let $0 < \epsilon, \ 0 < \lambda.$
\begin{itemize}
    \item We say that $x$ is $\epsilon$-feasible if $\max_i f_i(x) \le \epsilon$.  The quantity $\max_i f_i(x)$ is the \textit{infeasibility} of $x$.
    \item We say that $x$ is $\lambda$-superoptimal if $g(x) \ge \lambda \ + \ \max\left\{ g(x) \, : \, \, \text{s.t.} \, \eqref{eq:genpolyconst} \right\}$.
\end{itemize} 

\noindent {\bf Comment.} As defined, infeasibility and superoptimality are additive quantities. As such, they can be misleading if, for example, the coefficients defining the $f_i$ are very small (understating infeasibilities) or those defining $g$ are very large (overstating the superoptimality), or if the entries
in $x$ can be very large.  In the example discussed below and in 
Theorem \ref{infvsuper} all coefficients are small integers (or can be rescaled so that is the case) and the feasible region is
contained in a cube with sides of small integral magnitude. \\

\begin{figure}
    \centering
    \includegraphics[scale = 0.6]{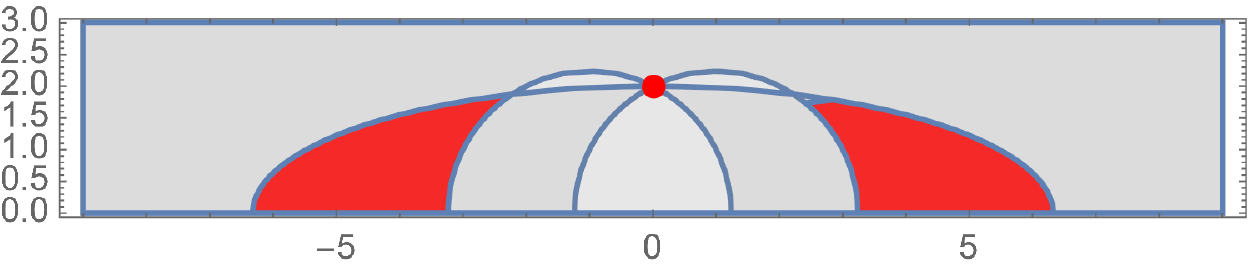}
    \includegraphics[scale = 0.6]{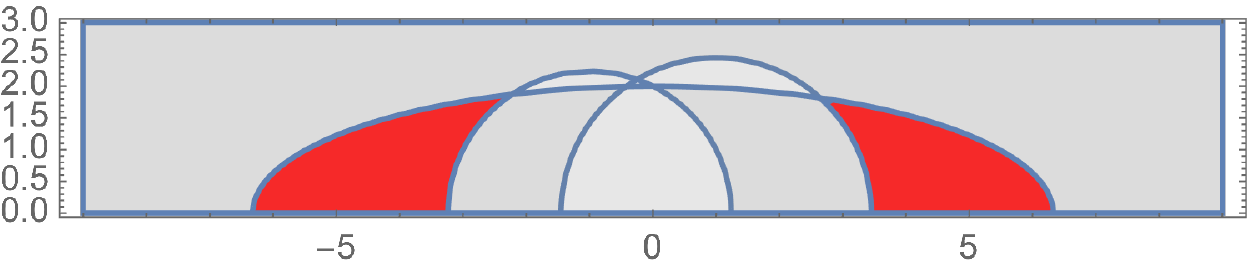}
    \caption{Red represents feasible region of system \eqref{eq:unlucky}.   Top: $\sigma = 0$.  Bottom: $\sigma = 1$.}
    \label{fig:quad-near-optimal}
\end{figure}

Now we turn to the construction. We first consider the system in variables $z_1, \, z_2$
\begin{subequations} \label{eq:unlucky}
\begin{align}
    & (z_1 - 1)^2 + z_2^2 \ \ge \ 5 + \sigma \label{eq:O1} \\
    & (z_1 + 1)^2 + z_2^2  \ \ge \ 5 \label{eq:O2} \\
    & \frac{z_1^2}{10} + z_2^2 \ \le \ 4 \label{eq:E} \\
    &  z_2 \ge 0 
\end{align}
\end{subequations}
where $0 \le \sigma \le 1$ is a parameter.  See Figure~\ref{fig:quad-near-optimal}.

\begin{lemma}\label{unlucky} (i) System \eqref{eq:unlucky} has no solutions 
with $0 < z_1 < 2$.  (ii) System \eqref{eq:unlucky} has no solutions with 
$-2 < z_1 < 0$.  (iii) If $\sigma = 0$ then $(z_1,z_2) = (0,2)$ is feasible.
(iv) If $\sigma > 0$ there is no feasible solution with $z_1 = 0$; further $z_2^2 \le 4 - 4/10$ in every feasible solution.
\end{lemma}
\begin{proof}  (i) Inequalities \eqref{eq:O1} and \eqref{eq:E} imply that  $9/10 z_1^2 - 2 z_1 - \sigma \ge 0$ which has no solutions with $0 < z_1 < 2$. (ii) Similar to (i). (iii) Clear by inspection. (iv) \eqref{eq:O1} and \eqref{eq:E} imply that $z_1 \neq 0$ when $\sigma > 0$. The rest follows from (i)-(ii) and \eqref{eq:E}. \qed
\end{proof}

Our complexity result concerns the following optimization problem which is constructed from a 3SAT instance $C_1 \wedge C_2 \ldots \wedge C_m$ using the recipe given by \eqref{eq:np-hard}; this is combined with a system \eqref{eq:unlucky}:
\begin{subequations}  \label{eq:combined}
  \begin{align}
 Z^*_2 := \ & \max \ z_2  \nonumber \\  
\text{s.t.} \quad  & (z_1 - 1)^2 + z_2^2 \ \ge \ 5 + s \label{eq:O1p} \\
    & (z_1 + 1)^2 + z_2^2  \ \ge \ 5 \label{eq:O2p} \\
    & \frac{z_1^2}{10} + z_2^2 \ \le \ 4 \label{eq:Ep} \\
    &  z_2 \ge 0 \label{eq:nonnegz2p}\\
    & \text{$(x_1, \ldots, x_{2n}, \gamma, \Delta, y_1, y_2, d_1, \ldots, d_n, s)$ feasible for \eqref{eq:np-hard}} \label{eq:np-hardp}
\end{align}
\end{subequations}
Constraints \eqref{eq:O1p}-\eqref{eq:nonnegz2p} are a copy of system \eqref{eq:unlucky} using the variable $s$ instead of $\sigma$; we stress that $s$ also appears in \eqref{eq:np-hardp}.\\

\noindent We now present a number of results regarding this system.  Let $0 < \epsilon < 1$ be given, and let $y^{(\epsilon)}$ be a rational approximation to the vector $y^*$ in Example \ref{example cubic v2} such that
(a) $h(y^{(\epsilon)}) < \epsilon$, (b) $y^{(\epsilon)} \in R_0$ and (c) $y^{(\epsilon)}$ has size polynomial in $\log \epsilon^{-1}$.  The fact that such a vector exists follows from the fact that $h$ is cubic a polynomial and therefore is Lipschitz-continuous on the bounded domain $R_0$.

The next statements focus on the vector $\bm{\omega^{(\epsilon)}}$ given by $x_j = 1 = -x_{n+j}$ for $j \in [n]$, $\Delta = 2$, $\gamma = 0$, 
$d_k = 0$ for $k \in n$, $y = y^{(\epsilon)}$, $s = 0$, $z_1 = 0$, and  $z_2 = 2$.\\
 
\noindent {\bf Statement 1.} The vector $\bm{\omega^{(\epsilon)}}$  is $\epsilon$-feasible for \eqref{eq:combined}. \\
\noindent{\em Proof.} Following the proof of Claim 1 in the proof of Theorem \ref{hard1}, it is clear that this vector satisfies constraints \eqref{eq:box}-\eqref{eq:gammadelta} as well as \eqref{eq:nastytiny}.  Moreover, the vector is also $\epsilon$-feasible for \eqref{eq:nasty} by construction of $y^{(\epsilon)}$.  Finally the vector is
feasible for \eqref{eq:O1p}-\eqref{eq:nonnegz2p} by Lemma \ref{unlucky} (ii). \qed 
\vspace{.1in}

\noindent {\bf Statement 2.} Suppose the formula $C_1 \wedge C_2 \ldots \wedge C_m$ is satisfiable. Then $Z^*_2 = 2$. As a result, $\bm{\omega^{\epsilon}}$    attains the optimal value for the corresponding instance of problem \eqref{eq:combined}.\\
\noindent {\em Proof.} By Claim 3, system  \eqref{eq:np-hard} has a feasible solution with $s = 0$.  Then Lemma \eqref{unlucky}(iii) completes the proof. \qed
\vspace{.1in}

\noindent {\bf Statement 3.} Suppose the formula $C_1 \wedge C_2 \ldots \wedge C_m$ is not satisfiable. Then $Z^*_2 \le \sqrt{4 - 4/10} < 1.9$.  As a result,
$\bm{\omega^{\epsilon}}$    is $0.2$-superoptimal for \eqref{eq:combined}.\\
\noindent {\em Proof.} Consider any feasible solution to \eqref{eq:combined}. By Claim 4 (a) $s > 0$.  Thus by Lemma \ref{unlucky}(iv) we have 
$z_2^2 \le 4 - 4/10$ which completes the proof. \qed
\vspace{.1in}
\noindent As a summary of the above discussion, we have:
\begin{theorem}\label{infvsuper}
Consider an instance of problem \eqref{eq:combined}.  Let $0 < \epsilon < 1$ and let $$\tilde v = (\tilde x_1, \ldots, \tilde x_{2n}, \tilde \gamma, \tilde \Delta, \tilde y_1, \tilde y_2, \tilde d_1, \ldots, \tilde d_n, \tilde s)$$
be $\epsilon$-feasible and of size polynomial in the size of the instance and $\log \epsilon^{-1}$.  It is strongly NP-hard to decide if $\tilde s \ge \frac{2}{1.9} Z^*_2$, even when it is known that either  
$\tilde s \ge \frac{2}{1.9} Z^*_2$ or $\tilde s = Z^*_2$.  
\end{theorem}

It is useful to contrast the situation described above with that arising in linear mixed-integer optimization.  Consider a standard (linear) MIP:
\begin{subequations} \label{eq:genmin}
\begin{align} 
  \textbf{MX:} \quad \max \ & \ c^T x \ + \ d^T y  \nonumber \\
    \text{s.t.} \ & \ A x + B y \le b  \nonumber\\
    \ & x \in \R^{n}, \ y \in \{0,1\}^m    \nonumber
\end{align}
\end{subequations}
Typically, solvers for MIP will produce vectors that are be slightly infeasible  due to round off errors in floating point arithmetic.  However, we can argue that the linear MIP case is not as problematic as the general nonlinear case (as exemplified in Theorem \ref{infvsuper}). To see this, suppose that $(x^*, y^*)$ is a candidate solution for {\bf MX} with $y^*$ binary\footnote{A numerical solver might yield a vector $y^*$ that is near-binary, i.e. within a small tolerance.  The analysis below is easily adjusted to handle such an eventuality.}, and suppose that this solution exhibits small infeasibilities -- could it be the case that this solution is far from any feasible solution for {\bf MX}, or that it has a large superoptimality relative to {\bf MX}?  

Let us provide a precise argument that shows that we can answer this question, in polynomial time. Consider the 
linear program  
\begin{subequations} \label{eq:genminfixed}
\begin{align} 
   \textbf{MXf}(y^*): \quad  \max \ & \ c^T x \nonumber\\
    \text{s.t.} \ & \ A x \le b - B y^* \nonumber\\
    \ & \ x \in \R^n,   \nonumber
\end{align}
\end{subequations}
\begin{itemize}
 \item [(a)] If {\bf MXf$(y^*)$} is \textit{infeasible}, then we can diagnose this fact \textit{in polynomial time}, and  conclude that the vector $(x^*, y^*)$ is `far' from any feasible solution for the {\bf MX}, namely its $L_\infty$-distance to any feasible point for {\bf MX} is at least $1$\footnote{Recall that $y^*$ is binary.} 
 \item[(b)] Suppose next that {\bf MXf$(y^*)$} is feasible. In this case we can show that the superoptimality of $(x^*, y^*)$ is upper bounded by a linear function of its infeasibility.  First, if $x^*$ is feasible for {\bf MXf$(y^*)$} then, certainly, $(x^*, y^*)$ is not superoptimal for {\bf MX}.  Suppose, on the other hand, that $x^*$  has infeasibility $\delta > 0$ for {\bf MXf$(y^*)$}.  Since {\bf MXf$(y^*)$} is feasible, its dual is bounded.  Let $\hat y \ge 0$ be an extreme point optimal solution to the dual of {\bf MXf$(y^*)$}.  Then
 $$ c^T x^* \ = \ \hat y^T A x^* \ \le \ \hat y^T(b - By^* + \delta \text{{\bf e}}) \ = \ \hat y^T(b - By^*) + \|\hat y\|_1 \delta$$
 where {\bf e} is the vector of 1s. In other words the superoptimality of $x^*$ for {\bf MXf$(y^*)$} --and thus, the superoptimality of $(x^*, y^*)$ for {\bf MX}-- is upper bounded by an expression of the form $\kappa(A,c) \delta$ where $\kappa(A,c)$ is a constant dependent on $A$ and $c$  which can furthermore be computed in polynomial time.  This upper bound on superoptimality amounts to a \textit{condition number} bound, and is completely missing from the general, nonlinear setting. 
 \end{itemize}

\subsection{Cubics and unbounded rays}
\label{sec:cubics-unbounded}

A \emph{ray} of a polyhedron $P \subseteq \R^n$ is a set of the form $R(x, d) := \{ x + \lambda d \mid \lambda \ge 0\}$, for some $x \in P$ and some nonzero $d$ in $\rec P$, where $\rec P$ denotes the recession cone of $P$ (see, e.g., \cite{Schrijver}). 
By definition of recession cone, every ray of a polyhedron $P$ is contained in $P$.
Next, we show that there can be a rational polyhedron $P$ and a cubic polynomial that is unbounded on some rays of $P$, but is bounded on all rational rays of $P$.
This is in contrast with the linear and quadratic case, where a rational unbounded ray is always guaranteed to exist. 
Formally, we say that a function $f\colon \R^n \to \R$ is \emph{bounded} on a set $P \subseteq \R^n$ if there exists $\phi \in \R$ such that $f(x) \le \phi$ for every $x \in P$.
Otherwise, for every $\phi \in \R$ there exists a vector $x \in P$ with $f(x) \ge \phi$; in other words $\sup_{x \in P} f(x) = + \infty$.
In the latter case we say that $f$ is \emph{unbounded} on $P$.

\begin{proposition}
\label{pr ray}
There exists a rational polyhedron $P \subseteq \R^3$ and a cubic polynomial $f \colon \R^3 \to \R$ that is unbounded on some rays of $P$, but for every ray $R(y, d)$ of $P$
such that $f$ is unbounded on $R(y, d)$, the vector $d$ is not rational.
\end{proposition}

\begin{proof}
We define the following three polynomial functions from $\R^3$ to $\R$
\begin{align*}
c(y) & := - 2 y_1^3 - y_2^3 + 6y_1y_2 y_3 - 4y_3^3, \\
q(y) & := y_1 y_3, \\
f(y) & := c(y) + q(y).
\end{align*}
The function $c(y)$ is a homogeneous cubic, while $q(y)$ is a homogeneous quadratic.
We also define the following two rational polyhedra
\begin{align*}
Q &:= \{y \in \R^3 \mid (y_1,y_2) \in [1.25,1.26] \times [1.58,1.59], \ y_3 = 1 \}, \\
P &:= \cone Q.
\end{align*}
Note that $- c(y_1,y_2,1)$ coincides with the cubic function considered in Example~\ref{example cubic v2}. 
Hence the maximum of $c(y)$ on $Q$ is zero and it is achieved only at the irrational vector $\tilde d := (2^{\frac 13}, 2^{\frac 23},1) \approx (1.2599,1.5874,1)$.

We first check that $f$ is unbounded on the ray of $P$ given by $R(0, \tilde d) = \{\lambda \tilde d \mid \lambda \ge 0\}$. 
We obtain 
\begin{align*}
c(\lambda \tilde d) & = \lambda^3 c(\tilde d) = 0, \\
q(\lambda \tilde d) & = \lambda^2 q(\tilde d), \\
f(\lambda \tilde d) & = c(\lambda \tilde d) + q(\lambda \tilde d) = \lambda^3 c(\tilde d) + \lambda^2 q(\tilde d) = \lambda^2 q(\tilde d).
\end{align*}
We observe that $f(\lambda \tilde d)$ is a quadratic univariate function in $\lambda$ and the leading coefficient is $q(\tilde d) = \tilde d_1 \tilde d_3 = 2^{\frac 13} \approx 1.2599 > 0$.
Thus $f \to +\infty$ along the ray $R(0, \tilde d)$.

Since the vector $\tilde d$ is irrational (and cannot be scaled to be rational), we now only need to consider rays of $P$ of the form $R(\bar y, \bar d)$,
for some $\bar y \in P$ and some $\bar d \in Q \setminus \{\tilde d\}$.
We fix $\bar y \in P$ and $\bar d \in Q \setminus \{\tilde d\}$ and evaluate the functions $c,q,f$ on $R(\bar y, \bar d)$:
\begin{align*}
c(\bar y + \lambda \bar d) & = \lambda^3 c(\bar d) + O(\lambda^2), \\
q(\bar y + \lambda \bar d) & = O(\lambda^2), \\
f(\bar y + \lambda \bar d) & = c(\bar y + \lambda \bar d) + q(\bar y + \lambda \bar d) = \lambda^3 c(\bar d) + O(\lambda^2).
\end{align*}
We observe that $f(\bar y + \lambda \bar d)$ is a cubic univariate function in $\lambda$ and the leading coefficient is $c(\bar d)$.
From Example~\ref{example cubic v2}, for every $\bar d \in Q \setminus \{\tilde d\}$ we have $c(\bar d) < 0$, thus $f \to -\infty$ along the ray $R(\bar y, \bar d)$.
\qed
\end{proof}

Note that the cubic function in the proof of Proposition~\ref{pr ray} grows only quadratically along the presented irrational unbounded ray.

\subsection{NP-Hard to determine if there is a rational unbounded ray}

The main result of this section is the following theorem.

\begin{theorem}\label{hard-unbounded} Consider unbounded optimization problems of the form 
$$
\max\{\pi(x) : \quad Ax \leq b\},
$$
where $\pi \in \Z[x_1, \dots, x_n]$ is a polynomial of degree three, and $A \in \Z^{m\times n}$, $b \in \Z^m$.   It is strongly NP-hard to test if there exists a rational ray on which the problem is unbounded.
\end{theorem}

The main reduction is from the problem 3SAT.  We use the same notation to describe an instance of 3SAT as above.   

Let $N \ge 1$ be an integer.  For now $N$ is generic; below we will discuss particular choices. Given an instance of 3SAT as above, we construct a system of quadratic inequalities on the following $2n + N + 7$ variables:
\begin{itemize}
\item For each literal $w_j$ we have a variable $x_j$; for $\bar w_j$ we use
  variable $x_{n+j}$.
  \item Additional variables $\gamma, \Delta, y_1, y_2, s,$ and $d_1, \ldots, d_N$.
\end{itemize}

Consider the polyhedral cone $K$ given by 
\begin{subequations}
\label{eq:cubic-unbounded-nphard}
\begin{align}
    -y_3 &\leq x_j \leq y_3 & j \in [2n]\\
    x_j + x_{n+j} &= 0 & j\in [n]\\
    0 &\leq y_3\\
    x_{i_1} + x_{i_2} + x_{i_3} &\geq -y_3 - \Delta & \forall C_i = (x_{i_1} \vee x_{i_2} \vee x_{i_3})\\
    0 &\leq \gamma \\
    0 & \leq \Delta \leq 2y_3\\
    \Delta + \frac{\gamma}{2} &\leq 2y_3\\
    1.259 y_3 - \gamma &\leq y_1 \leq 1.26 y_3\\
    1.587 y_3 & \leq y_2 \leq 1.59 y_3.
\end{align}
\end{subequations}

Now consider the optimization problem
\begin{equation}
\label{eq:unbounded-opt}
    \max -n^6 y_3^3 + n^5 y_3 \sum_{j=1}^n x_j^2 + f(y) \ \ s.t. \ \ (y,x,\Delta, \gamma) \in K.
\end{equation}
Let $\pi(y,x,\Delta, \gamma)$ be the objective function.  Let $\tilde d = \begin{pmatrix} 2^{1/3} \\ 2^{2/3} \\ 1\end{pmatrix}$.

\begin{theorem}\label{unbounded-hard1} Let $n \ge 3$.  Consider an instance \eqref{eq:3sat}  of 3-SAT and the corresponding  optimization problem ~\eqref{eq:unbounded-opt} with feasibility region $K$ in \eqref{eq:cubic-unbounded-nphard}. 
\begin{itemize}
    \item [\textbf{(a)}] The objective function 
    $\pi\to +\infty$ along the ray given by $y = \tilde d y_3$, $x_j = y_3$, $x_{n+j} = 0$ for all $j \in [n]$, $\Delta = 2y_3$, $\gamma = 0$ where we let $y_3 \in [0,+\infty)$.

    \item [\textbf{(b)}] 
    Suppose formula \eqref{eq:3sat} is satisfiable.  Then there is a rational, feasible ray (i.e., contained in $K$) 
    over which $\pi \to +\infty$.  The ray has of bit encoding size $4n + \mathcal{L}$, where $\mathcal{L}$ is a fixed constant independent of $n$.
    
    \item [\textbf{(c)}] Suppose formula \eqref{eq:3sat} is not satisfiable. 
    If $\pi \to +\infty$ along some ray, then this ray is irrational.
\end{itemize}
\end{theorem}

\begin{proof} \hspace{.1in} \\ 
  \noindent \textbf{(a)}  
  This statement is clear.
  
  \hspace{.1in} \\ 
  \noindent \textbf{(b)} 
  Define the ray as follows.  First, set $\bar \Delta = 0$, $\bar \gamma = 4$, and for each $j \in [n]$, let $\bar x_j = \pm 1$ if literal $j$ is true or false and $\bar x_{n+j} = 0$.  
Finally, let $\bar y_1$,$\bar y_2$ be rational such that $1.259 - 4 \leq \bar y_1 \leq 1.26$ and $1.587 \leq \bar y_2 \leq 1.59$ and $-2\bar y_1^3 - \bar y_2^3 + 6 \bar y_1 \bar y_2 - 4 > 0$.  Such values $(\bar y_1 , \bar y_2)$ exist (see Observation~\ref{ob1}).
Then along the ray given by 
$
\begin{pmatrix}
y\\x\\\Delta\\\gamma
\end{pmatrix} 
= 
\begin{pmatrix}
\bar y\\ \bar x\\ \bar \Delta\\ \bar \gamma
\end{pmatrix} y_3,
$
we have $\pi = f(y) \to + \infty$ (cubically).
  
  \hspace{.1in} \\ 
  \noindent \textbf{(c)} 
  First, $y_3 \to +\infty$ along this ray (else every variable is bounded).  So we can write this ray as 
$$
\alpha^0 + \begin{pmatrix} \hat y\\ \hat x \\ \hat \Delta \\ \hat \gamma 
\end{pmatrix}
y_3
$$
where $\hat y_3 = 1$ and for some $\alpha^0 \in \R^{2n+5}$.
Since $\hat \Delta \geq 0$ and $\hat \Delta +  \tfrac{\hat \gamma}{2} \leq 2$, we have $\hat \gamma \leq 4$.  Hence, along this ray, $f(y) \leq 12 y_3^3 + \mathcal{O}(y_3^2)$.   This follows from Observation~\ref{ob1}. 

It follows that $-n^6 + n^5 \sum_{j=1}^n \hat x_j^2 \geq -12$, else $\pi \to -\infty$ cubically in $y_3$ along this ray.  Thus $\hat x_j^2 \geq 1 - \mathcal{O}(\tfrac{1}{n^5}) \forall j$.

Consider the truth assignment where literal $j \in [n]$ is true if and only if $\hat x_j > 0$.  So there exists a clause $i$ such that $\hat x_{i_1} + \hat x_{i_2} + \hat x_{i_3} \leq -3 + \mathcal{O}(\tfrac{1}{n^5})$, and thus $\hat \Delta \geq 2 - \mathcal{O}(\tfrac{1}{n^5})$.  and so $\hat \gamma \leq \mathcal{O}(\tfrac{1}{n^5})$.  In that case, 

\begin{equation}
    \begin{cases}
    1.259 y_3 - \gamma \leq y_1 \leq 1.26 y_3\\
    1.587 y_3 \leq y_2 \leq 1.59 y_3
    \end{cases}
    \Rightarrow 
        \begin{cases}
    1.259 y_3 - \mathcal{O}(\tfrac{1}{n^5}) \leq y_1 \leq 1.26 y_3\\
    1.587 y_3 \leq y_2 \leq 1.59 y_3,
    \end{cases}
\end{equation}
and so $f(y) \to -\infty$ cubically in $y_3$ unless $\hat y = \begin{pmatrix}
2^{1/3} \\ 2^{2/3} \\1
\end{pmatrix}$
and also $-n^6 y_3^3 + n^5 y_3 \sum_{j=1}^n x_j^2$ either goes to $-\infty$ or is equal to $0$.  So we must have $\hat y = \begin{pmatrix}
2^{1/3} \\ 2^{2/3} \\1
\end{pmatrix}$.
\hfill \qed
  \end{proof}
  The proof of Theorem~\ref{hard-unbounded} follows directly from the above result.

\section{Short certificate of feasibility: an almost feasible point}
\label{sec:NP-fixed-nd}

In this section we are interested in the existence of \emph{short} certificates of feasibility for systems of polynomial inequalities, i.e., certificates of feasibility of size bounded by a polynomial in the size  of the system.

First, for completeness, we present a result on Lipschitz continuity of a polynomail on a box. 
\begin{lemma}[Lipschitz continuity of a polynomial on a box]
\label{lem:lipschitz}
Let $g \in \Z[x_1, \dots, x_n]$ be a polynomial of degree at most $d$ with coefficients of absolute value at most $H$.
Let $y,z \in [-M,M]^n$ for some $M > 0$.
Then
\begin{equation}
|g(y) - g(z)|\leq L \|y - z\|_\infty
\end{equation}
where $L := nd H M^{d-1} (n+d)^{d-1}$.
\end{lemma}

\begin{proof}
Let $g(x) = \sum_{I \in \N^n, \normone I \le d} c_I x^I$.
By the fundamental theorem of calculus,
\begin{equation*}
g(z) = g(y) + \int_{\lambda =0}^{\lambda = 1} \nabla g(y + \lambda (z-y))^\top (z-y)\, d\lambda.
\end{equation*}
Therefore, 
an upper bound on $|g(y) - g(z)|$ can be obtained by bounding the quantity $|\nabla g(y + \lambda(z-y))^\top (z-y)|$ with $\lambda \in [0,1]$.
In the remainder of the proof we derive such a bound. 

Note that, for every $i=1,\dots,n$, we can write 
\begin{align}
    \label{eq nabla}
    (\nabla g(x) )_i = \sum_{I \in \N^n, \normone I \le d-1} \tilde c_{I,i} x^I,
\end{align}
where $|\tilde c_{I,i}| \leq d H$.
Hence, for any 
$I \in \N^n$ with $\normone I \le d-1$ and $x \in [-M,M]^n$,
we have
$
|\tilde c_{I,i} x^I| \leq d H M^{d-1}.
$

Next we bound the number of terms in the summation in \eqref{eq nabla}.
A \emph{weak composition} of an integer $q$ into $p$ parts is a sequence of $p$ non-negative integers that sum up to $q$.
Two sequences that differ in the order of their terms define different weak compositions.
It is well-known that the number of weak compositions of a number $q$ into $p$ parts is $\binom{q+p-1}{p-1} = \binom{q+p-1}{q}$.
(For more details on weak compositions see, for example, \cite{HeuManBook}.)
We obtain that the number of terms in the summation in \eqref{eq nabla} is bounded by
$$
\sum_{i=0}^{d-1} \binom{i+n-1}{i}
\le
\sum_{i=0}^{d-1} \binom{d+n-2}{i}
\le
(d+n-1)^{d-1}
\le
(n+d)^{d-1},
$$
where in the second inequality we used the binomial theorem.

We obtain that for every $x \in [-M,M]^n$,
$$
\|\nabla g(x)\|_\infty \leq d H M^{d-1} (n+d)^{d-1}.
$$
Therefore, for any 
$\lambda \in [0,1]$,
\begin{align*}
|\nabla g(y + \lambda(z-y))^\top (z-y)| 
& \leq 
n \cdot
\|\nabla g(y + \lambda(z-y))\|_\infty \cdot \|y-z\|_\infty \\
& \leq 
n \cdot
d H M^{d-1} (n+d)^{d-1}
\cdot \|y-z\|_\infty.
\end{align*}
\qed
\end{proof}

\begin{figure}
    \centering
    a.)\  \includegraphics[scale=0.4]{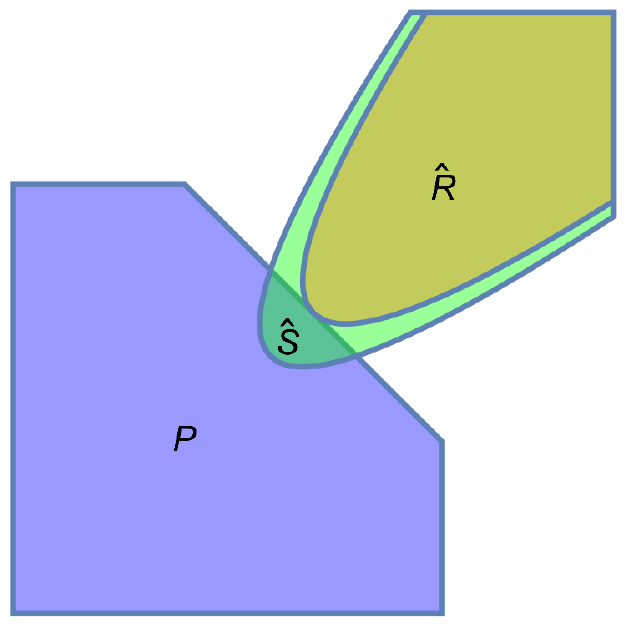} \hspace{2cm} b.) \  \includegraphics[scale=0.4]{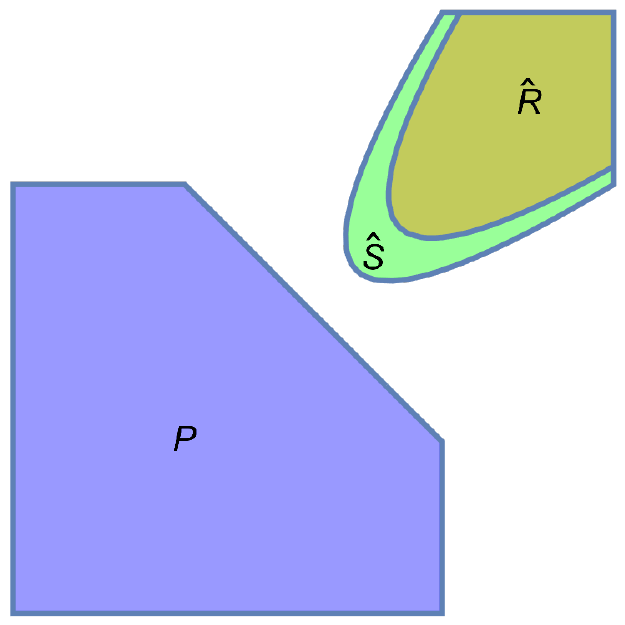}
    \caption{Structure of sets in Proposition~\ref{prop relaxed vector} where we define sets $\hat R := \{x \in \R^n \mid g_j(x) \le 0, j \in [\ell] \}, \
\hat S := \{x \in \R^n \mid \ell \delta g_j(x) \le 1, j \in [\ell] \}$.  Hence $R = P \cap \hat R$ and $S = P \cap \hat S$. a.) Example where $R$ is feasible.  b.) Example where $R$ is infeasible.   }
    \label{fig: propo}
\end{figure}

Our first result relies on a boundedness assumption, which seems to be common in results of this type. It is not clear if that assumption can be removed.
We refer the reader to Figure~\ref{fig: propo} for an illustration of the sets considered in the next result.

\begin{proposition}
\label{prop relaxed vector}
Let $f_i \in \Z[x_1, \dots, x_n]$, for $i \in [m]$, of degree one.
Let $g_j \in \Z[x_1, \dots, x_n]$, for $j \in [\ell]$, of degree bounded by an integer $d$.
Assume that the absolute value of the coefficients of $f_i$, $i \in [m]$, and of $g_j$, $j \in [\ell]$, is at most $H$.
Let $\delta$ be a positive integer.
Let $P := \{ x \in \R^n \mid f_i(x) \le 0, i \in [m]\}$ and 
consider the sets
\begin{align*}
R := \{x \in P \mid g_j(x) \le 0, j \in [\ell] \}, \
S := \{x \in P \mid \ell \delta g_j(x) \le 1, j \in [\ell] \}.
\end{align*}
Assume that $P$ is bounded.
If $R$ is nonempty, then there exists a rational vector in $S$ of size bounded by a polynomial in $n,d, \log \ell, \log H, \log \delta$.
\end{proposition}

\begin{proof}
Since $P$ is bounded, it follows from Lemma 8.2 in \cite{BerTsi97} that $P \subseteq [-M,M]^n$, where $M = (nH)^n$.
Let $L$ be defined as in Lemma~\ref{lem:lipschitz}, i.e.,
\begin{align*}
L := n d H M^{d-1} (n+d)^{d-1} = n d H (nH)^{n(d-1)} (n+d)^{d-1}.
\end{align*}
Note that $\log L$ is bounded by a polynomial in $n,d,\log H$.
Let $\varphi := \ceil{L M \ell \delta}$.
Therefore $\log \varphi$ is bounded by a polynomial in $n,d, \log \ell, \log H, \log \delta$.

We define the following $(2\varphi)^n$ boxes in $\R^n$ with  $ j_1,\dots,j_n \in \{-\varphi,\dots, \varphi-1\}$:
\begin{align}
\label{eq boxes}
\C_{j_1,\dots,j_n} := 
\left\{x \in \R^n \mid \frac{M}{\varphi} j_i \le x_i \le \frac{M}{\varphi} (j_i + 1), i \in [n] \right\}.
\end{align}
Note that the union of these $(2\varphi)^n$ boxes is the polytope $[-M,M]^n$ which contains the polytope $P$.
Furthermore, each of the $2n$ inequalities defining a box \eqref{eq boxes} has size polynomial in $n,d, \log \ell, \log H, \log \delta$.

Let $\tilde x$ be a vector in $R$, which exists because $R$ is assumed nonempty.
Since $\tilde x \in P$, there exists a box among \eqref{eq boxes}, say $\tilde \C$, that contains $\tilde x$.
Let $\bar x$ be a vertex of the polytope $P \cap \tilde \C$. 
Since each inequality defining $P$ or $\tilde \C$ has size polynomial in $n,d, \log \ell, \log H, \log \delta$, it follows from Theorem 10.2 in \cite{Schrijver} that also $\bar x$ has size polynomial in $n,d, \log \ell, \log H, \log \delta$.

To conclude the proof of the theorem we only need to show $\bar x \in S$.
Since $\bar x, \tilde x \in \tilde \C$, we have $\|\bar x - \tilde x\|_\infty \leq \frac M \varphi$.
Then, from Lemma~\ref{lem:lipschitz}
we obtain that for each $j \in [\ell]$,
$$
| g_j(\bar x) - g_j(\tilde x)| \le L \|\bar x - \tilde x\|_\infty \le \frac {LM} {\varphi} 
\le \frac{1}{\ell \delta}.
$$
If $g_j(\bar x) \le 0$ we directly obtain $g_j(\bar x) \le \frac{1}{\ell \delta}$ since $\frac{1}{\ell \delta} > 0$.
Otherwise we have $g_j(\bar x) > 0$.
Since $g_j(\tilde x) \le 0$, we obtain $g_j(\bar x) \le | g_j(\bar x) - g_j(\tilde x)| \le \frac{1}{\ell \delta}$.
We have shown that $\bar x \in S$, and this concludes the proof of the theorem.
\qed
\end{proof}

We will be using several times the functions $\epsilon$ and $\delta$ defined 
as follows:
\begin{align*}
\epsilon(n,m,d,H) 
& := (2^{4-\frac{n}{2}} \max\{H, 2n + 2m\} d^n)^{-n 2^n d^n}, \\
\delta(n,m,d,H) 
& := \ceil{2\epsilon^{-1}(n,m,d,H)}
= \ceil{2(2^{4-\frac{n}{2}} \max\{H, 2n + 2m\} d^n)^{n 2^n d^n}}.
\end{align*}

A fundamental ingredient in our arguments is
the following result by Geronimo, Perrucci, and Tsigaridas, which follows from Theorem 1 in \cite{Geronimo2013}.

\begin{theorem}
\label{thm:epsilon}
Let $n \ge 2$.
Let $g,f_i \in \Z[x_1, \dots, x_n]$, for $i \in [m]$, of degree bounded by an even integer $d$.
Assume that the absolute value of the coefficients of $g,f_i$, $i \in [m]$, is at most $H$.
Let 
$
T := \{x \in \R^n \mid f_i(x) \le 0, i \in [m] \},
$
and let $C$ be a compact connected component of $T$.
Then, the minimum value that $g$ takes over $C$,
is either zero, or its absolute value is greater than or equal to $\epsilon(n,m,d,H)$.
\end{theorem}

Using Theorem~\ref{thm:epsilon} we obtain the following lemma.
\begin{lemma}
\label{lem one ine}
Let $n \ge 2$.
Let $g,f_i \in \Z[x_1, \dots, x_n]$, for $i \in [m]$, of degree bounded by an even integer $d$.
Assume that the absolute value of the coefficients of $g,f_i$, $i \in [m]$, is at most $H$.
Let $\delta := \delta(n,m,d,H).$
Let $T:= \{x \mid f_i(x) \le 0, i \in [m] \}$ and consider the sets
\begin{align*}
R  := \{x \in T \mid g(x) \le 0\}, \
S  := \{x \in T \mid \delta g(x) \le 1 \}.
\end{align*}
Assume that $T$ is bounded.
Then $R$ is nonempty if and only if $S$ is nonempty.
\end{lemma}

\begin{proof}
Since $\delta > 0$ we have $R \subseteq S$, therefore if $R$ is nonempty also $S$ is nonempty. Hence we assume that $S$ is nonempty and we show that $R$ is nonempty.

Since $S$ is nonempty, there exists a vector $\bar x \in T$ with $g(\bar x) \le 1/\delta < \epsilon(n,m,d,H)$.
Let $C$ be a connected component of $T$ containing $\bar x$.
Since $T$ is compact, we have that $C$ is compact as well.
In particular, the minimum value that $g$ takes over $C$ is less than $\epsilon(n,m,d,H)$.
The contrapositive of Theorem~\ref{thm:epsilon} implies that the minimum value that $g$ takes over $C$ is less than or equal to zero.
Thus there exists $\tilde x \in C$ with $g(\tilde x) \le 0$.
Hence the set $R$ is nonempty.
\qed
\end{proof}

\ifdefined\journalversion

From Lemma~\ref{lem one ine} we obtain the following result.
\begin{proposition}
\label{prop many ine}
Let $n \ge 2$.
Let $f_i,g_j \in \Z[x_1, \dots, x_n]$, for $i \in [m]$, $j \in [\ell]$, of degree bounded by an even integer $d$.
Assume that the absolute value of the coefficients of $f_i$, $i \in [m]$, and of $g_j$, $j \in [\ell]$, is at most $H$.
Let $\delta := \delta(n,m+\ell,2d,\ell H^2).$
Let $T := \{x  \mid f_i(x) \le 0, i \in [m] \}$ and consider the sets
\begin{align*}
R := \{x \in T \mid g_j(x) \le 0, j \in [\ell] \}, \
S  := \{x \in T \mid \ell \delta g_j(x) \le 1, j \in [\ell]  \}.
\end{align*}
Assume that $T$ is bounded.
Then $R$ is nonempty if and only if $S$ is nonempty.
\end{proposition}

\begin{proof}
Since $\ell \delta > 0$ we have $R \subseteq S$, therefore if $R$ is nonempty also $S$ is nonempty. Hence we assume that $S$ is nonempty and we show that $R$ is nonempty.

Let $\bar x \in S$, and define the index set $J := \{j \in [\ell] : g_j(\bar x) > 0\}.$
We introduce the polynomial function $g \in \Z[x_1, \dots, x_n]$ defined by
$g(x) := \sum_{j \in J} g_j^2(x).$
Note that the degree of $g$ is bounded by $2d$.
The absolute value of the coefficients of each $g^2_j$ is at most $H^2$, hence the absolute value of the coefficients of $g$ is at most $\ell H^2$.
Next, let $T':= \{x \in T \mid \, g_j(x) \le 0, j \in [\ell]\setminus J\}$ and 
\begin{align*}
R'  := \{x \in T' \mid g(x) \le 0 \}, \ \ 
S'  := \{x \in T' \mid \delta g(x) \le 1 \}.
\end{align*}

First, we show that the vector $\bar x$ is in the set $S'$, implying that $S'$ is nonempty.
Clearly $\bar x \in T$, and for every $j \in [\ell]\setminus J$ we have that $g_j(\bar x) \le 0$, thus we have $\bar x \in T'$.
For every $j \in J$, we have $0< g_j(\bar x) \le \frac{1}{\ell \delta}$, and since  $\ell \delta \ge 1$, we have $0 < g_j^2(\bar x) \le \frac{1}{\ell \delta}$.
Thus, we obtain $g(\bar x) \le \frac{\ell}{\ell \delta} = \frac{1}{\delta}$.
We have thus proved $\bar x \in S'$, and so $S'$ is nonempty.

Next, we show that the set $R'$ is nonempty.
To do so, we apply Lemma~\ref{lem one ine} to the sets $T',R',S'$.
The number of inequalities that define $T'$ is a number $m'$ with $m \le m' \le m + \ell$.
The degree of $f_i,g_j,g$, for $i \in [m]$, $j \in [\ell]\setminus J$, is bounded by $2d$.
The absolute value of the coefficients of $f_i,g_j,g$, for $i \in [m]$, $j \in [\ell]\setminus J$, is at most $\ell H^2$.
Since the function $\delta(n,m,d,H)$ 
is increasing in $m$ and $m' \le m + \ell$, we obtain from Lemma~\ref{lem one ine} that $R'$ is nonempty if and only if $S'$ is nonempty.
Since $S'$ is nonempty, we obtain that $R'$ is nonempty.

Finally, we show that the set $R$ is nonempty.
Since $R'$ is nonempty, let $\tilde x \in R'$.
From the definition of $R'$ we then know $\tilde x \in T$, $g_j(\tilde x) \le 0$, for $j \in [\ell]\setminus J$, and $g(\tilde x) \le 0$.
Since $g$ is a sum of squares, $g(\tilde x) \le 0$ implies $g(\tilde x) = 0$, and this in turn implies $g_j(\tilde x) = 0$ for every $j \in J$.
Hence $\tilde x \in R$, and $R$ is nonempty.
\qed
\end{proof}

Proposition~\ref{prop many ine} and Proposition~\ref{prop relaxed vector} directly yield our following main result.

\begin{theorem}[Certificate of polynomial size]
\label{th short certificate}
Let $n \ge 2$.
Let $f_i \in \Z[x_1, \dots, x_n]$, for $i \in [m]$, of degree one.
Let $g_j \in \Z[x_1, \dots, x_n]$, for $j \in [\ell]$, of degree bounded by an even integer $d$.
Assume that the absolute value of the coefficients of $f_i$, $i \in [m]$, and of $g_j$, $j \in [\ell]$, is at most $H$.
Let $\delta := \delta(n,m+\ell,2d,\ell H^2)$.
Let $P := \{ x \in \R^n \mid f_i(x) \le 0, i \in [m]\}$ and consider the sets
\begin{align*}
R  := \{x \in P \mid g_j(x) \le 0, j \in [\ell] \}, \
S  := \{x \in P \mid \ell \delta g_j(x) \le 1, j \in [\ell]. \end{align*}
Assume that $P$ is bounded.
Denote by $s$ the maximum number of terms of $g_j$, $j \in [\ell]$, with nonzero coefficients.
If $R$ is nonempty, then there exists a rational vector in $S$ of size bounded by a polynomial in $d,\log m, \log \ell, \log H$, for $n$ fixed.
This vector is a certificate of feasibility for $R$ that can be checked in a number of operations that is bounded by a polynomial in $s, d, m, \ell, \log H$, for $n$ fixed.
\end{theorem}

\begin{proof}
From the definition of $\delta$,
we have that $\log \delta$ is bounded by a polynomial in $d,\log m, \log \ell, \log H$, for $n$ fixed.
From Proposition~\ref{prop relaxed vector}, there exists a vector $\bar x \in S$ of size bounded by a polynomial in $d,\log m, \log \ell, \log H$, for $n$ fixed.
Such a vector is our certificate of feasibility.
In fact, from Proposition~\ref{prop many ine} (applied to the sets $T = P, R,S$), we know that $S$ nonempty implies $R$ nonempty.

To conclude the proof we  bound the 
number of operations 
needed to 
check if the vector $\bar x$ is in $S$ by
substituting $\bar x$ in the $m+\ell$ inequalities defining $S$.

The absolute value of the coefficients of $f_i(x) \le 0$, $i \in [m]$, is at most $H$.
Thus, it can be checked that $\bar x$ satisfies these $m$ inequalities in a number of operations that is bounded by a polynomial in $d, m, \log \ell, \log H$, for $n$ fixed.

Next, we focus on the inequalities $\ell \delta g_j(x) \le 1$, $j \in [\ell]$.
Note that the total number of terms of $\ell \delta g_j$, $j \in [\ell]$, with nonzero coefficients is bounded by $s \ell$.
The logarithm of the absolute value of each nonzero coefficient of $\ell \delta g_j(x)$, $j \in [\ell]$, is bounded by $\log(\ell \delta H)$, which in turn is bounded by a polynomial in $d,\log m, \log \ell, \log H$, for $n$ fixed.
Therefore, it can be checked that $\bar x$ satisfies these $\ell$ inequalities in a number of operations that is bounded by a polynomial in $s, d, \ell, \log m, \log H$, for $n$ fixed.
\qed
\end{proof}

In particular, Theorem~\ref{th short certificate} implies that polynomial optimization is in NP, provided that we fix the number of variables.As mentioned in the introduction, this fact is not new. 
In fact, it follows from Theorem~1.1 in Renegar~\cite{Ren92f} that the problem of deciding whether the set $R$, as defined in Theorem~\ref{th short certificate}, is nonempty can be solved in a number of operations that is bounded by a polynomial in $s, d, m, \ell, \log H$, for $n$ fixed.
Therefore, Renegar's algorithm, together with its proof, provides a certificate of feasibility of size bounded by a polynomial in the size of the system, which in turns implies that the decision problem is in NP.

The main advantages of Theorem~\ref{th short certificate} over Renegar's result are that 
(i) our certificate of feasibility is simply a vector in $S$ of polynomial size, and 
(ii) the feasibility of the system can be checked by simply plugging the vector into the system of inequalities defining $S$.
The advantages of Renegar's result over our Theorem~\ref{th short certificate} are:
(iii) Renegar does not need to assume that the feasible region is bounded, while we do need that assumption, and
(iv) Renegar shows that the decision problem is in P, while we show that it is in the larger class NP.
Note that in ~\cite{renegar92}, Renegar shows how to produce points that are within some distance from a feasible solution.

\section{Existence of rational feasible solutions}
\label{sec:rational-feasible}
In this section we present a set of results that will be used to argue that rational solutions exist to certain feasibility problems.  We begin with a number of auxiliary results that we will rely on.

The first was known by \href{https://scholarship.richmond.edu/cgi/viewcontent.cgi?article=1113&context=masters-theses}{Nicol\`o Fontana Tartaglia}.  
It shows that in a univariate cubic polynomial, 
shifting by a constant allows us to assume that 
the $x^2$ term has a zero coefficient.
\begin{lemma}[Rational shift of cubic]
\label{lem:rational-shift}
Let $f(x) = ax^3 + bx^2 + cx + d$.  Then $f(y-\tfrac{b}{3a}) = a y^3 + \tilde c y + \tilde  d$ where $\tilde c = \frac{ 27 a^2 c-9 a b^2}{27 a^2}$ and $\tilde d = \frac{27 a^2 d-9 a b c+2 b^3}{27 a^2}$.
\end{lemma}
\begin{proof}
Consider an assignment of variables $x = y - s$ where $y$ is a new variable and $s$ is the shift.
\begin{align}
    f(y-s) &= a(y-s)^3 + b(y-s)^2 + c(y-s) + d\\
    &= a y^3+ (b-3 a s)y^2+  \left(3 a s^2-2 b s+c\right)y-a s^3+b s^2-c s+d
\end{align}
Hence, setting $s = \frac{b}{3a}$, we have
\begin{equation}
    f\left(y - \frac{b}{3a}\right) = a
   y^3 +\frac{ 27 a^2 c-9 a b^2}{27 a^2}y + \frac{27 a^2 d-9 a b c+2 b^3}{27 a^2}.
\end{equation}
\qed
\end{proof}

The next lemma provides bounds on the roots of a univariate polynomial.  We attribute this result to Cauchy; a proof can be found in Theorem 10.2 of~\cite{Basu2006}.

\begin{lemma}[Cauchy - size of roots]
\label{lem Cauchy}
Let $f(x) = a_n x^n + \dots + a_1 x + a_0$, where $a_n, a_0 \neq 0$.  Let $\bar x \neq 0$ such that $f(\bar x) = 0$.  Then  $\tfrac{1}{\delta} \leq |\bar x| \leq M$, where
$$
M = 1 + \max\left\{\left|\tfrac{a_0}{a_n}\right|, \dots, \left|\tfrac{a_{n-1}}{a_n}\right|\right\}, \ \ \delta = 1 + \max\left\{\left|\tfrac{a_1}{a_0}\right|, \dots, \left|\tfrac{a_{n}}{a_0}\right|\right\}.
$$
Furthermore, if $a_i \in \Z$ for $i=0,1,\dots, n$, then $\lceil M \rceil$ and $\lceil \delta \rceil$ are integers of size polynomial in the sizes of $a_0, a_1, , \dots, a_n$.
\end{lemma}
The next lemma is a  special case of Theorem~2.9 in \cite{Alb11}. 
\begin{lemma}
\label{lem:rational-radical-soln}
Let $n \geq 1$.
Let $r_i \in \Q_+$ for $i=0,1, \dots, n$, $q_i \in \Q_+$ for $i=1, \dots, n$.  If 
$
    \sum_{i=1}^n r_i \sqrt{q_i} = r_0,
$
then $\sqrt{q_i} \in \Q$ for all $i=1, \dots, n$.  
Furthermore, the size of $\sqrt{q_i}$ is polynomial in the size of $q_i$.
\end{lemma}

Next we show that local minimizers of separable cubic polynomials are rational provided that the function value is rational.

\begin{theorem}[Rational local minimum]
\label{thm:rational-roots}
Let $f(x) = \sum_{i=1}^n f_i(x_i)$ where $f_i(x_i) = a_i x_i^3 + b_ix_i^2 +c_ix_i + d_i \in \Z[x_i]$ and $a_i \neq 0$ for all $i \in [n]$.
Assume that the absolute value of the coefficients of $f$ is at most $H$.
Suppose $x^*$ is the unique local minimum of $f$ and $\gamma^* := f(x^*)$ is rational.
Then $x^*$ is rational and has size that is  
polynomial in $\log H$ and in the size of $\gamma^*$.
\end{theorem}
\begin{proof}
For every $i \in [n]$, Let $\tilde c_i,\tilde d_i$ be defined as in Lemma~\ref{lem:rational-shift}, let $g_i(y_i) := a_i y_i^3 + \tilde c_i y_i$, and define $g(y) = \sum_{i=1}^n g_i(y_i)$.
Then $y^* \in \R^n$ defined by $y_i^* := x_i^* + \tfrac{b_i}{3a_i}$, $i \in [n]$, is the unique local minimum of $g(y)$ and $g(y^*) = \gamma^* - \sum_{i=1}^n \tilde d_i$.

We now work with the gradient.
Since $y^*$ is a local minimum of $g$, we have $\nabla g(y^*) = 0$.
Since $g(y)$ is separable, we obtain that for every $i \in [n]$,
\begin{equation}
\label{eq this one and no other}
     g_i'(y^*_i) = 0 \ \ \Rightarrow \ \ y^*_i = \pm \sqrt{\frac{-\tilde c_i}{3a_i}}.
\end{equation}
Furthermore, we will need to look at the second derivative.  
Since $y^*$ is a local minimizer, then $\nabla^2g(y^*) \ge 0$.  
Again, since $g(y)$ is separable, this implies that $g''_i(y_i) \geq 0$ for every $i \in [n]$.  
Hence we have
\begin{equation}
    g_i''(y_i) \geq 0 \ \ \Rightarrow \ \ 6a_i y^*_i \geq 0 \ \ \Rightarrow \ \ a_i \left( \pm \sqrt{\frac{-\tilde c_i}{3a_i}} \right) \geq 0.
\end{equation}
Also, notice that we must have $\frac{-\tilde c_i}{3a_i} \geq 0$ for  $\sqrt{\frac{-\tilde c_i}{3a_i}}$ to be a real number.
Thus, 
\begin{equation}
\label{eq:signs}
    \textrm{sign}(- \tilde c_i) = \textrm{sign}(a_i) = \textrm{sign}\left( \pm \sqrt{\frac{-\tilde c_i}{3a_i}}\right).
\end{equation}
Finally, we relate this to $g(y^*)$.
\begin{align}
    \gamma^* - \sum_{i=1}^n \tilde d_i = g(y^*) 
    = \sum_{i=1}^n \left(a_i (\tfrac{-\tilde c_i}{3a_i}) \left(\pm \sqrt{\tfrac{-\tilde c_i}{3a_i}}\right) + \tilde c_i\left(\pm \sqrt{\tfrac{-\tilde c_i}{3a_i}}\right)\right) 
    = -\tfrac{2}{3} \sum_{i=1}^n |\tilde c_i| \sqrt{\frac{-\tilde c_i}{3a_i}},
\end{align}
where the last equality comes from comparing the signs of the data from \eqref{eq:signs}.
Hence, we have 
\begin{equation}
\label{eq calc}
\sum_{i=1}^n |\tilde c_i| \sqrt{\frac{-\tilde c_i}{3a_i}} = -\tfrac{3}{2} (\gamma^* - \sum_{i=1}^n \tilde d_i).
\end{equation}
By Lemma~\ref{lem:rational-radical-soln},  for every $i \in [n]$, $\sqrt{\frac{-\tilde c_i}{3a_i}}$ is rational  and has size polynomial in $\log H$ and in the size of $\gamma^*$.
From \eqref{eq this one and no other}, so does $y^*$, and hence $x^*$.
\qed
\end{proof}

We are now ready to prove our main result of this section.

\begin{theorem}
\label{th onetwo}
Let $n \in \{1,2\}$.
Let $f_i \in \Z[x_1,\dots,x_n]$, for $i \in [m]$, of degree one.
Let $g(x) = \sum_{i=1}^n (a_i x_i^3 + b_i x_i^2 +c_i x_i + d_i) \in \Z[x_1,\dots,x_n]$ with $a_i \neq 0$ for $i \in [n]$.
Assume that the absolute value of the coefficients of $g,f_i$, $i \in [m]$, is at most $H$.
Consider the set
\begin{align*}
R & := \{x \in \R^n \mid g(x) \le 0, \ f_i(x) \le 0, i \in [m] \}.
\end{align*}
If $R$ is nonempty, then it contains a rational vector of size bounded by a polynomial in $\log H$.
This vector provides a certificate of feasibility for $R$ that can be checked in a number of operations that is bounded by a polynomial in $m, \log H$.
\end{theorem}
 
\begin{proof}
Define $P := \{ x \in \R^n \mid f_i(x) \le 0, i \in [m]\},$ let $x^*$ be a vector in $P$ minimizing $g(x)$, and let $\gamma^* := g(x^*)$.

Since $R$ is nonempty, we have $\gamma^* \le 0$.
If $\gamma^* = 0$, then Theorem~\ref{thm:rational-roots} implies that the size of $x^*$ is bounded by a polynomial in $\log H$, thus the result holds.
Therefore, in the remainder of the proof we assume $\gamma^* < 0$.

Without loss of generality, we assume that that $x^*$ is in the interior of $P$.  Otherwise, if $x^*$ is contained in a lower dimensional face of $P$, we can project into that face.
In particular, we assume that $x^*$ is the unique local minimum of $g$.

Note that if $x^*$ is in a 0-dimensional face of $P$, then since $P$ is a rational polyhedron, $x^*$ rational and of size bounded by a polynomial in $\log H$.

\textbf{Claim:} \emph{There exits an integer $\delta$ of polynomial size bounded by a polynomial in $\log H$ such that $|\gamma^*| < \tfrac{1}{\delta}$.}
We prove separately the cases $n=1,2$.

\emph{Claim proof.} 
Case $(n=1)$.  Let $\tilde c,\tilde d$ be defined as in Lemma~\ref{lem:rational-shift}.
Following the calculation of Theorem~\ref{thm:rational-roots}, \eqref{eq calc}
$
\gamma^* - \tilde d = - \tfrac{2}{3} |\tilde c| \left(\pm \sqrt{\tfrac{-\tilde c}{3a}}\right).
$
Hence, 
$
    (\gamma^* - \tilde d)^2 = -\frac{4 \tilde c^3}{27a},
$
that is, $\gamma^*$ is a non-zero root of the above quadratic equation.
From Lemma~\ref{lem Cauchy}, we have that $|\gamma^*| \ge \frac{1}{\delta}$, where $\delta$ is an integer and $\log \delta$ is bounded by a polynomial in $\log H$.

Case $(n=2)$.  For every $i \in [n]$, let $\tilde c_i,\tilde d_i$ be defined as in Lemma~\ref{lem:rational-shift}.
Following the calculation of Theorem~\ref{thm:rational-roots}, \eqref{eq calc}
$
\sum_{i=1}^2 |\tilde c_i| \sqrt{\frac{-\tilde c_i}{3a_i}} 
= - \tfrac{3}{2} (\gamma^* - \sum_{i=1}^2 \tilde d_i).
$
Squaring both sides we obtain
$$
\sum_{i=1}^2 \frac{-\tilde c_i^3}{3a_i} 
+ 2 |\tilde c_1| |\tilde c_2| \sqrt{\frac{\tilde c_1 \tilde c_2}{ 9 a_1 a_2}}
= \tfrac{9}{4} (\gamma^* - \sum_{i=1}^2 \tilde d_i)^2.
$$
If we isolate the square root, and then square again both sides of the equation, we obtain that $\gamma^*$ is a non-zero root of a quartic equation with rational coefficients.
From Lemma~\ref{lem Cauchy}, we have that $|\gamma^*| \ge \frac{1}{\delta}$, where $\delta$ is an integer and $\log \delta$ is bounded by a polynomial in $\log H$.
This concludes the proof of the claim. \hfill $\diamond$

Since $\gamma^* <0$, we have thereby shown that $x^*$ is a vector in $P$ satisfying $g(x^*) \le - \frac{1}{\delta}$.

Since $x^*$ is on the interior of $P$, it must satisfy $\nabla g(x^*) = 0$.  Hence, for $i\in[n]$, $x^*_i$ is a root of the quadratic equation 
$$
(a_i x_i^3 + b_i x_i^2 +c_i x_i + d_i)' = 3 a_i x_i^2 + 2 b_i x_i +c_i = 0,
$$
thus using again Lemma~\ref{lem Cauchy} we obtain that $-M \le x_i^* \le M$, where $M$ is an integer and $\log M$ is bounded by a polynomial in $\log H$.
We apply Proposition~\ref{prop relaxed vector} to the polytope $\{ x \in \R^2 \mid f_i(x) \le 0, i \in [m], \ -M \le x_i \le M, i \in [n]\}$, with $\ell:=1$, and with $g_1(x) := g(x) + \frac{1}{\delta}$.
Proposition~\ref{prop relaxed vector} then implies that there exists a vector $\bar x \in P$ with $g_1(\bar x) \le \frac{1}{\delta}$, or equivalently $g(\bar x) \le 0$, of size bounded by a polynomial in $\log H$.
Such a vector is our certificate of feasibility.

\smallskip

To conclude the proof for both cases $n=1$ and $n=2$, we bound the number of operations needed to check if $\bar x$ is in $R$ by substituting $\bar x$ in the $m+1$ inequalities defining $R$.
It is simple to check that $\bar x$ satisfies these $m+1$ inequalities in a number of operations that is bounded by a polynomial in $m, \log H$.
\qed
\end{proof}

Example~\ref{example cubic v2}  shows that Theorem~\ref{th onetwo} cannot be generalized to non-separable bi-variate cubics.  Furthermore, that example can be easily be lifted to a separable cubic in three dimensions.  To see this, let $y_3 = y_2 - y_1$.  Then 
$$
3y_3^2 - 3y_1^2 - 3 y_2^2 = - 6y_1y_2. 
$$
Thus, we can replace the mixed term with separable terms, provided we lift this to three dimensions and add an extra equation.

\section{Unbounded rays for cubic objectives}
In this section we present conditions on the existence of rays of a polyhedron, along which a polynomial function is unbounded. As in the remainder of the paper, special attention is given to the rationality of these rays.
In this section we denote by $\norm{\cdot}$ the 2-norm, although any norm would suffice.

We will use the following standard lemma.
For a proof (and a more general statement) we refer the reader to \cite{AndBelShi82}, \cite[Lemma 2.1]{Kla19}.
\begin{lemma}
\label{lem Kla}
Let the set $P \subset \R^n$ be nonempty and closed, let $f \colon P \to \R$ be continuous on $P$ and suppose that $f$ is unbounded on $P$.
Then there exists a sequence $\{x^j\}_{j \in \N}$ of vectors in $P$ such that 
\begin{itemize}
\item[(i)] $f(x^{j-1}) < f(x^j) \to +\infty$, 
\item[(ii)] $\norm{x^{j-1}} < \norm{x^j} \to + \infty$, and
\item[(iii)] $\forall x \in P$, $\norm{x} < \norm{x^j} \Rightarrow f(x) < f(x^j)$.
\end{itemize}
\end{lemma}

We are now ready to prove the first result of this section, which discusses unbounded functions on polyhedra and on its rays.
The proof technique is essentially the one used in
classic proofs of Frank-Wolfe type theorems in cubic optimization, see \cite{AndBelShi82}, \cite{BelAnd93}, \cite{Kla19}, and most closely follows that in~\cite{Kla19}. 

\begin{theorem}
\label{th ray}
Let $P \subseteq \R^n$ be a 
polyhedron and let $f \colon \R^n \to \R$ be a polynomial of degree at most three.
If $f$ is unbounded on $P$, then there exists 
a ray of $P$ over which $f$ is unbounded.
\end{theorem}

\begin{proof}
The proof is by contradiction.
Thus we assume that there exists a counterexample to the theorem, which consists of a polyhedron $P \subseteq \R^n$ and a polynomial $f \colon \R^n \to \R$ of degree at most three such that $f$ is unbounded on $P$ but is bounded over each ray $R(x, d)$ of $P$.
Among all counterexamples, we consider one where the polyhedron $P$ has minimal dimension. Clearly, we have $\dim P \ge 1$.
Using an affine function that maps the affine hull of $P$ onto $\R^{\dim P}$, we can assume without loss of generality that $n = \dim P$, i.e., that $P$ is full-dimensional.

By Lemma~\ref{lem Kla}, there exists a sequence $\{x^j\}_{j \in \N}$ of vectors in $P$ such that 
\begin{align}
& f(x^{j-1}) < f(x^j) \to +\infty, \label{eq lem i} \\
& \norm{x^{j-1}} < \norm{x^j} \to + \infty, \label{eq lem ii} \\
& \forall x \in P, \quad \norm{x} < \norm{x^j} \Rightarrow f(x) < f(x^j). \label{eq three}
\end{align}
For every $j \in \N$ we define the scalar $\tau_j := \norm{x^j}$ and the direction vector $d^j := x^j / \tau_j \in \R^n$.
We can then rewrite \eqref{eq lem ii} in the form
\begin{align}
\label{eq tau}
0 < \tau_{j-1} < \tau_j \to +\infty. 
\end{align}
Clearly we have $\norm{d^j} = 1$, thus the vectors $d^j$ lie on the unit sphere, which is a compact set.
The Bolzano-Weierstrass Theorem implies that the sequence $\{d^j\}_{j \in \N}$ has a convergent subsequence whose limit is in the unit sphere. We denote by $d$ this limit, and from now on we only consider without loss of generality such a subsequence, thus we can write 
$d^j \to d$.

Next we show $d \in \rec P.$
Let $Ax \le b$ be a system of linear inequalities defining $P$, i.e., $P = \{x \mid Ax \le b\}$.
The definition of $d^j$ and 
$x^j \in P$ imply that $Ad^j = Ax^j /\tau_j \le b / \tau_j$.
By taking the limits and using $d^j \to d$ and \eqref{eq tau}, we obtain $Ad \le 0$, i.e., $d \in \rec P$.

We now show that there exists an index $s_1 \in \N$ such that
\begin{align}
\label{eq four}
\forall j \ge s_1, \quad x^j \in \intr P,
\end{align}
where $\intr P$ denotes the interior of $P$.
In order to prove this, it suffices to show that only finitely many vectors $x^j$ are in $P \setminus \intr P$.
We prove this latter statement by contradiction, and so we assume that infinitely many vectors $x^j$ are in $P \setminus \intr P$.
We note that the set $P \setminus \intr P$ is the union of the finitely many faces $F$ of $P$ with $\dim F < n$.
Hence, there exists a face $F$ of $P$ with $\dim F < n$ that contains infinitely many vectors $x^j$, which implies that the function $f$ is unbounded on $F$.
The minimality of our counterexample implies that the theorem is true for the polyhedron $F$ and the polynomial $f$, and so there exists a ray of $F$ over which $f$ is unbounded. 
This is a contradiction because each ray of $F$ is also a ray of $P$.
Hence, \eqref{eq four} is shown.

Next, we show that that there exists an index $s_2 \in \N$ such that
\begin{align}
\label{eq six}
\forall j \ge s_2 \text{ and } \forall \mu \in (0,1], \quad \norm{x^j -\mu d}<\norm{x^j}.
\end{align}
To prove \eqref{eq six}, we first observe that there exists an index $s_2 \in \N$ such that 
\begin{align*}
\forall j \ge s_2, \quad \tau_j>1  \text{ and }  \norm{d^j}-\norm{d^j-d} > 0.
\end{align*}
This follows from $\tau_j \to +\infty$ in \eqref{eq tau}, 
$d^j \to d$ 
and $\norm{d^j} = 1$.
Using the triangle inequality we then obtain 
\begin{align*}
\norm{x^j - \mu d} & = \norm{\tau_jd^j - \mu d} 
\le \norm{\tau_j d^j - \mu d^j} + \norm{\mu d^j - \mu d} = \\
& = \tau_j \norm{d^j} - \mu \norm{d^j} + \mu \norm{d^j - d}
< \tau_j \norm{d^j}
= \norm{x^j}.
\end{align*}
This completes the proof of \eqref{eq six}.

We now set $s := \max\{s_1,s_2\}$ and prove the following two properties:
\begin{align}
\label{eq nine}
\exists \mu_s  \in (0, 1] : & \quad f (x^s - \mu_s d) < f (x^s), \\
\label{eq ten}
\exists \lambda_s > 0 : & \quad f (x^s + \lambda_s d) < f (x^s).
\end{align}
We start by proving \eqref{eq nine}.
From \eqref{eq four}, $x^s \in \intr P$, which implies that there exists 
$\mu_s  \in (0, 1]$ such that $x^s - \mu_s d \in P$.
From \eqref{eq six} we have $\norm{x^s -\mu_s d}<\norm{x^s}$ and so \eqref{eq three} implies \eqref{eq nine}.
Next we show \eqref{eq ten}.
Since $x^s \in P$ and $d \in \rec P$, we can consider the ray $R(x^s,d)$
of $P$. 
On the ray $R(x^s,d)$, the polynomial $f(x)$ is not constant, due to \eqref{eq nine}, and it is  bounded by assumption.
Hence $\lim_{\lambda \to +\infty} f(x^s + \lambda d) = -\infty$ which implies \eqref{eq ten}.

Now we define for each $j > s$ the vector $v^j := (x^j - x^s) / \norm{x^j - x^s}$
and the restriction $F_j(\lambda)\colon \R \to \R$ of $f(x)$ to the line
$\{x^s + \lambda v^j \mid \lambda \in \R \}$, i.e.,
$F_j(\lambda) := f(x^s + \lambda v^j)$.
Next we show that $v^j \to d$.
We can write $v^j$ as follows
\begin{align*}
v^j & = \frac{x^j - x^s}{\norm{x^j - x^s}} = \frac{x^j}{\norm{x^j - x^s}} - \frac{x^s}{\norm{x^j - x^s}} = \\
& 
= \frac{\tau_j d^j}{\norm{\tau_j d^j - \tau_s d^s}} - \frac{x^s}{\norm{x^j - x^s}} 
= \frac{d^j}{\norm{d^j - (\tau_s/\tau_j)d^s}} - \frac{x^s}{\norm{x^j - x^s}}.
\end{align*}
The first addend goes to $d$ since $d^j \to d$, $\norm{d^j} = 1$, and $\tau_j \to +\infty$ from \eqref{eq tau}.
The second addend goes to zero since $\norm{x^j} \to + \infty$ in \eqref{eq lem ii} implies $\norm{x^j - x^s} \to +\infty$ for $j \to \infty$.

Next, we show that there exists $s' > s$ such that, for all $j \ge s'$,
\begin{align}
\label{eq eleven a}
& F_j(-\mu_s) < F_j(0), \\
\label{eq eleven b}
& F_j(\lambda_s) < F_j(0), \\
\label{eq thirteen}
\exists \lambda_j > \lambda_s : \quad & F_j(\lambda_j) > F_j(0).
\end{align}
To prove inequalities \eqref{eq eleven a} and \eqref{eq eleven b}, we only need to notice that since $f$ is continuous and $v^j \to d$, there exists $s' > s$ such that for all $j \ge s'$, $d$ in \eqref{eq nine} and \eqref{eq ten} can be replaced by $v^j$:
\begin{align*}
& F_j(-\mu_s) = f(x^s - \mu_s v^j) < f (x^s) = F_j(0), \\
& F_j(\lambda_s) = f (x^s + \lambda_s v^j) < f (x^s) = F_j(0).
\end{align*}
Without loss of generality we can assume that $s'$ is large enough such that, for every $j \ge s'$, we have $\lambda_s < \norm{x^j - x^s}$.
This is because $\lambda_s$ is fixed and $\norm{x^j - x^s} \to +\infty$ for $j \to \infty$. 
We let $\lambda_j := \norm{x^j -x^s}$ and obtain $\lambda_j > \lambda_s$.
Inequality \eqref{eq thirteen} is then obtained by applying \eqref{eq lem i} and $j > s$:
\begin{align*}
\exists \lambda_j > \lambda_s : \quad & F_j(\lambda_j) = f(x^j) > f (x^s) = F_j(0).
\end{align*}

Thus, for $j\ge s'$, \eqref{eq eleven a}, \eqref{eq eleven b}, and \eqref{eq thirteen}
imply that $F_j(\lambda)$ increases somewhere in the interval $[-\mu_s , 0]$, then decreases somewhere in $[0, \lambda_s]$, and again increases somewhere in $[\lambda_s, \lambda_j]$. 
In particular, $F_j(\lambda)$ cannot be a linear or quadratic function.
Since $f(x)$ is a cubic polynomial, we obtain that $F_j(\lambda)$ is a genuinely cubic polynomial.
Hence $F_j(\lambda)$ must increase on the interval $[\lambda_j, +\infty)$, i.e.,
$F_j(\lambda)$ is unbounded on the half-line $L_j := \{x^s + \lambda(x^j - x^s) \mid \lambda \ge 0\}$.
Since by assumption $f$ is bounded over each ray of $P$,
then
each half-line $L_j$, for $j \ge s',$ has to leave the polyhedron $P$ at some vector $y^j \in P \setminus \intr P$ with $f(y^j)>f(x^j)$.

Now consider the sequence of vectors $\{y^j\}_{j \ge s'}$. 
Then, there exists a face $F$ of $P$ with $\dim F < n$ that contains infinitely many vectors $y^j$ with $j \ge s'$.
This in particular implies that $f$ is unbounded on $F$.
The minimality of our counterexample implies that the theorem is true for the polyhedron $F$ and the polynomial $f$, and so there exists a ray of $F$ over which $f$ is unbounded. 
This is a contradiction because each ray of $F$ is also a ray of $P$.
\qed
\end{proof}

It should be noted that analogues of Theorem~\ref{th ray} are known to hold for linear and quadratic functions \cite{Schrijver,Vav90,Pia2016}. 

Next, we discuss the tightness of Theorem~\ref{th ray} with respect to the degree of the polynomial function. 
Namely, we show that Theorem~\ref{th ray} does not hold if the polynomial $f$ of degree at most three is replaced with a quartic polynomial.
The example given in the next proposition is inspired by a similar example by Frank and Wolfe \cite{FraWol56} (see also \cite{Kla19}).

\begin{proposition}
\label{th no ray}
There exists a quartic polynomial $f \colon \R^2 \to \R$ unbounded on $\R^2$ such that $f$ is bounded over every ray of $\R^2$.
\end{proposition}

\begin{proof}
We define the following quartic polynomial function from $\R^2$ to $\R$
\begin{align*}
f(y_1,y_2) & := y_2 - (y_2 - y_1^2)^2.
\end{align*}
To check that $f$ is unbounded on $\R^2$, it suffices to consider the vectors that satisfy $y_2 = y_1^2$ as $y_1 \to + \infty$.

Hence, to prove the proposition we only need to check that $f$ is bounded over every ray of $\R^2$.
Hence, we let $\bar y \in \R^2$, $\bar d \in \R^2 \setminus \{0\}$ and we evaluate the function $f$ on the ray $R(\bar y, \bar d)$:
\begin{align*}
f(\bar y + \lambda \bar d) 
& = \bar y_2 + \lambda \bar d_2 - (\bar y_2 + \lambda \bar d_2 - (\bar y_1 + \lambda \bar d_1)^2)^2 \\
& = \bar y_2 + \lambda \bar d_2 - (- \lambda^2 \bar d_1^2 + \lambda (\bar d_2 - 2 \bar y_1 \bar d_1) + \bar y_2 - \bar y_1^2)^2 \\
& = - \lambda^4 \bar d_1^4 
+ 2 \lambda^3 \bar d_1^2 (\bar d_2 - 2 \bar y_1 \bar d_1) 
+ \lambda^2 (2 \bar d_1^2 (\bar y_2 - \bar y_1^2) - (\bar d_2 - 2 \bar y_1 \bar d_1)^2) + O(\lambda).
\end{align*}
We observe that $f(\bar y + \lambda \bar d)$ is a quartic univariate function in $\lambda$.
If $\bar d_1 \neq 0$, then the leading term is $- \lambda^4 \bar d_1^4$.
In this case, since $\bar d_1^4 > 0$, we obtain $f \to -\infty$ along the ray $R(\bar y, \bar d)$.
We now consider the remaining case $\bar d_1 = 0$.
In this case the leading term is $\lambda^2 (2 \bar d_1^2 (\bar y_2 - \bar y_1^2) - (\bar d_2 - 2 \bar y_1 \bar d_1)^2) = - \lambda^2 \bar d_2^2$.
Since $\bar d$ is nonzero and $\bar d_1 = 0$, we obtain $\bar d_2 \neq 0$.
Thus $\bar d_2^2 > 0$, and $f \to -\infty$ along the ray $R(\bar y, \bar d)$.
\qed
\end{proof}

Proposition~\ref{pr ray} shows that in Theorem~\ref{th ray} there might not exist any \emph{rational} ray along which the cubic is unbounded, even if we further assume that P is rational.
In particular, one might wonder if it is possible to construct an example similar to the one given in the proof of Proposition~\ref{pr ray}, but where the cubic function grows cubically along an unbounded ray.
The next proposition implies that this is not possible.
In fact, in Proposition~\ref{prop rat cub unb ray} we show that, if there exists an unbounded ray along which a cubic grows cubically, then there exists also a \emph{rational} unbounded ray.

\begin{proposition}[Rational cubically unbounded ray]
\label{prop rat cub unb ray}
Let $P$ be a rational polyhedron and let $\bar x \in P$ and $\bar v \in \rec{P}$.
Suppose that $f(\bar x + \lambda \bar v) \to +\infty$ as $\lambda \to +\infty$ with $f(\bar x + \lambda \bar v) = \Theta(\lambda^3)$.  Then, for any $\epsilon > 0$ there exist $\tilde x, \tilde v \in \Q^n$ with $\tilde x \in P$, $\tilde v \in \rec{P}$ such that $\norm{\bar x - \tilde x} < \epsilon$, $\norm{\bar v - \tilde v} < \epsilon$ and $f(\tilde x + \lambda \tilde v) \to +\infty$ as $\lambda \to +\infty$ with $f(\tilde x + \lambda \tilde v) = \Theta(\lambda^3)$.
\end{proposition}

\begin{proof}
Suppose that $\bar v$ is not rational. Since $\bar v$ is not rational, but $P$ is a rational polyhedron, then $\dim(\rec{P}) > 0$. Let $\Delta \in \R^n$ and consider
$$
f(\bar x + \lambda (\bar v + \Delta)) = \lambda^3f_3(\Delta) + \lambda^2f_2(\Delta) + \lambda f_1(\Delta) + f_0(\Delta)
$$
where $f_i \colon \R^n \to \R$ is of degree $i$.  By assumption, $f_3(0) > 0$.  By continuity of $f_3$, there exists a $\delta > 0$ such that $f_3(\Delta) > 0$ for any $\norm{\Delta}< \epsilon$.  
Since $\rec{P} \cap (\R^n \setminus \Q^n)$ is dense in $\rec{P} \cap \R^n$, there exists a (many) perturbations $\Delta$ from $\bar v$ such that $\tilde v := \bar v + \Delta \in \Q^n$ and $\norm{\Delta} < \epsilon$.
Since $f_i(\Delta)  > 0$, then $f(\bar x + \lambda \tilde v) \to +\infty$ as $\lambda \to +\infty$.

A similar argument applies to $\bar x$.
\qed
\end{proof}

We remark that, if the dimension is considered fixed, then in polynomial time we can determine if there exists a ray that is unbounded cubically.  This follows from Renegar's work~\cite{renegar1998}, bounds on minimum values of polynomial~\cite{Geronimo2013}, and optimizing the cubic term of the objective function over the boundary of a ball intersected with the recession cone of the feasible region.

To conclude this section, we observe that every cubic function has a direction in $\R^n$ that grows cubically.
\begin{observation}[All cubics have a direction that grows cubically]
Suppose $f \colon \R^n \to \R$ is a degree-3 polynomial.  Then there exists $\bar x, \bar v \in \R^n$ such that $f(\bar x + \lambda \bar v) = \Theta(\lambda^3)$.
\end{observation}
\begin{proof}
Let $f = f_0 + f_1 + f_2 + f_3$ where $f_i$ is a homogeneous polynomial of degree-$i$.  

Since $f$ is degree 3, there exists some $\bar v \in \R^n$ such that $f_3(\bar v) \neq 0$.  Since $f_3$ is cubic, it is symmetric about the origin, so we can assume that $f_3(\bar v) > 0$. Since $f_3$ is homogeneous of degree 3, $f_3(\lambda \bar v) = \lambda^3 f_3(\bar v)$. 

Hence, $f(\lambda \bar v) = \lambda^3 f_3(\bar v) + O(\lambda^2)$.  The lemma follows with $\bar x = 0$. \hfill \qed
\end{proof}

\section{Examples}
We conclude with a few related examples of irrationality and exponential size of solutions that can arise in seemingly innocent types of sets.  
\begin{example}[System of convex quadratics all of whose feasible solutions are large]
\label{example n quadratics}
Consider the system of inequalities 
$y_1 \ge 2$, $y_{i+1} - y_i^2 \ge 0$ for $i =1,\dots,n-1$.
Then each feasible vector satisfies $y_n \ge 2^{2^{n-1}}$.  Attributed to Ramana \cite{ramana} and Khachiyan, see \cite{letchfordparkes}, \cite{alizadeh}.\hfill $\diamond$
\end{example}

\begin{example}[A bounded feasible region QCQP whose solution requires exponentially many bits]
  \label{badboy}
\end{example}
\noindent Consider the optimization problem
\begin{subequations}  \label{bby}
  \begin{align}
\max & \ \quad x_2  \nonumber \\  
\text{s.t.} & \ (x_1 - 1)^2  + \ x_2^2 \ - \ d_N^2 \ \ge \ 3 \label{bby1}\\
&  \  (x_1 + 1)^2  + \ x_2^2 \ \ge \ 3 \label{bby2} \\
&  \ \frac{x_1^2 }{10} + x_2^2 \ \le \ 2 \label{bby3} \\
&  \ d_1 + d_N  = \frac{1}{2}, \quad 0 \le d_1, \quad d_{i}^2 \le d_{i+1} \ (1 \le i \le N-1) \label{sneaky}
\end{align}
\end{subequations}

Suppose we allow for $\epsilon$-feasible solutions, in the additive sense.  We
will show that unless $\epsilon < 2^{-2^{N}}$, there is an $\epsilon$-feasible solution to problem \eqref{bby} that attains value $\sqrt{2}$, whereas the true value of the problem is less than $1.23$.

We begin with the latter statement.  First, it is clear that \eqref{sneaky} implies that $d_N > 0$.  Armed with this fact, we will show that \eqref{bby1}-\eqref{bby3} imply that $x_2 < 1.3$. To see this, note that
\eqref{bby2} and \eqref{bby3} together imply that 
$$ 3 \le (x_1 + 1)^2 + 2 - \frac{x_1^2 }{10},$$
or $ 0 \le x_1(\frac{9}{10}x_1 + 2)$. Hence, if $x_1 < 0$ then $x_1 \le -20/9$
and so by \eqref{bby3}, $x_2 < \sqrt{2 - (20/9)^2/10} \approx 1.228$. Likewise, \eqref{bby1} and \eqref{bby3} together imply that
$$ 3 + d_N^2 \le (x_1 - 1)^2 + 2 - \frac{x_1^2 }{10},$$
or $d_N^2 \le x_1(\frac{9}{10}x_1 - 2)$.  Hence, if $x_1 \ge 0$, then $x_1 > 0$ (because $d_N > 0$) and $x_1 \ge 20/9$, and as above by \eqref{bby3}, $x_2 < \sqrt{2 - (20/9)^2/10} \approx 1.228$. Thus, in any case, $x_2 \le 1.228$.

At the same time, the vector given by $x_1 = 0$, $x_2 = \sqrt{2}$, $d_1 = \frac{1}{2}$, $d_i = \frac{1}{2^{2^i}}$ ($2 \le i \le N-1$) and $d_N = 0$ satisfies (exactly) all constraints \eqref{bby1}-\eqref{sneaky}, except for $d_{N-1}^2 \le d_N$, which it violates by $2^{-2^N}$.  This concludes the proof.\hfill $\diamond$

\begin{example}[An SOCP all of whose feasible solutions are irrational]
  \label{ex:irrational}
\end{example}

\noindent 
Let $(a,b,c,d)$ be a Pythagorean quadruple, i.e., $a^2 + b^2 + c^2 = d^2$, and all are integers.  Consider the system
\begin{subequations}
  \begin{align}
      & \sqrt{x_1^2 + x_2^2} \ \le \ x_0  \\
    & \sqrt{x_0^2 + x_3^2} \ \le \ d  \\
    & a \le x_1, \, b \le x_2, \, c \le x_3, 
  \end{align}
\end{subequations}
In any feasible solution we have 
$$
a^2 + b^2 \leq x_1^2 + x_2^2 \leq \, x_0^2 \, \leq d^2 - x_3^2 \leq d^2 - c^2 = a^2 + b^2.
$$
Hence, in any feasible solution, $x_0 = \sqrt{a^2 + b^2}$.  Choosing $(a,b,c,d) = (1,2,2,3)$, we have $x_0 = \sqrt{5}$.

Likely, examples of this type are already known.\hfill $\diamond$

%
%
%
\bibliographystyle{splncs04.bst}
\bibliography{references}
\end{document}